\documentclass[english, 11pt]{article}
\usepackage{amsmath}
\usepackage{amsthm}

\makeatletter

\usepackage{fullpage,graphicx,psfrag,amsfonts,verbatim, url}

\usepackage{abstract} 

\usepackage{parskip}
\begingroup
\makeatletter
   \@for\theoremstyle:=definition,remark,plain\do{\expandafter\g@addto@macro\csname th@\theoremstyle\endcsname{\addtolength\thm@preskip\parskip
     }}
\endgroup
\usepackage{thmtools}

\declaretheoremstyle[
headfont = \sffamily\bfseries,
notefont=\normalfont, notebraces={(}{)},
bodyfont=\normalfont\itshape,
headformat=\NAME~\NUMBER\NOTE
]{plain}
\declaretheoremstyle[
headfont = \sffamily\bfseries,
notefont=\normalfont, notebraces={(}{)},
bodyfont=\normalfont\itshape,
headpunct={},
headformat=\NAME~\NUMBER\NOTE
]{notnum}

\declaretheorem[numberwithin=section,style=plain]{definition}
\declaretheorem[numberwithin=section,style=plain]{assumption}

\usepackage{titlesec}
\def\sectionfont{\sffamily\Large\bfseries\boldmath}
\def\subsectionfont{\sffamily\large\bfseries\boldmath}
\def\paragraphfont{\sffamily\normalsize\bfseries\boldmath}
\titleformat*{\section}{\sectionfont}
\titleformat*{\subsection}{\subsectionfont}
\titleformat*{\subsubsection}{\paragraphfont}
\titleformat*{\paragraph}{\paragraphfont}
\titleformat*{\subparagraph}{\paragraphfont}
\usepackage[small,labelfont={bf,sf}]{caption}  

\usepackage{multirow}
\usepackage{enumitem}
\setlist{nolistsep}

\usepackage{amsmath, amssymb}
\usepackage{graphicx}
\usepackage{booktabs}
\usepackage{bm, bbm}
\usepackage{amsfonts}
\usepackage{xcolor}
\usepackage[margin=1in]{geometry}
\usepackage[round]{natbib}
\usepackage[hidelinks]{hyperref}
\usepackage{subcaption}
\usepackage[export]{adjustbox}
\usepackage{placeins}
\usepackage{longtable}
\usepackage{mathtools}
\usepackage[usenames,dvipsnames]{pstricks}
\usepackage{epsfig}
\usepackage{pst-grad} \usepackage{pst-plot}

\usepackage{array}
\newcolumntype{C}[1]{>{\centering\arraybackslash}p{#1}}

\allowdisplaybreaks

\makeatother

\usepackage{babel}

\begin{document}

\title{\textsf{\textbf{A Data-Driven Linear Programming Model for Energy-Optimal Metro Timetables}}}
\author{
  \small{\textbf{Qiuyu Fan}}\\
  \small{Department of Computational Applied Mathematics and Operations Research,}\\
  \small{Rice University, Houston, TX, USA} \\
  \small{\texttt{qf3@rice.edu}} \\
  \and
  \small{\textbf{J. Kevin Tobin}}\\
  \small{Hitachi Rail Canada, Toronto, ON, Canada} \\
  \small{\texttt{kevin.tobin@hitachirail.com}} \\
  \and
  \small{\textbf{Bart P.G. Van Parys}} \\
  \small{Centrum Wiskunde \& Informatica (CWI), Amsterdam, The Netherlands} \\
  \small{\texttt{bart.van.parys@cwi.nl}} \\
  \and
    \small{\textbf{Shuvomoy Das Gupta}}\\
  \small{Department of Computational Applied Mathematics and Operations Research,}\\
  \small{Rice University, Houston, TX, USA} \\
  \small{\texttt{sd158@rice.edu}} \\
}

\date{}

\maketitle
\begin{abstract}
We propose a data-driven linear programming model to compute energy-optimal
timetables for networks operating under communications-based train control
(CBTC); the model was developed in collaboration with Hitachi Rail Canada,
the largest provider of CBTC systems worldwide. Our model minimizes the
network's \emph{effective energy consumption}, defined as the total traction energy consumed
by the trains minus the regenerative braking energy transferred from braking
trains to nearby accelerating trains. Our main modeling contribution is to show that the signed overlap between a braking and an accelerating phase is concave in the event times, so maximizing the transferred energy is a convex problem that can be reformulated into a linear program via hypograph constraints. In contrast, prior work either captures this synchronization structure with binary variables or other nonconvex constraints, or relies on multiple stages
interleaved with simulation. Our model computes an energy-optimal timetable subject to the operational constraints of the railway network in under a second on a laptop for every instance studied and predicts the effective energy consumption of the network without requiring
time-consuming simulations. We demonstrate the effectiveness of our model on two real-world CBTC networks. First, on 11 full-day operational instances of Shanghai Metro
Line 8 with 1,000--1,332 active trains, the model computes energy-optimal
timetables in under 0.7 seconds per instance and predicts 19.5\%--28.2\%
reductions in effective energy consumption relative to the existing
operational timetables; the industrial physics-based simulator
\textsf{SPSIM} independently corroborates these reductions. Second, on 17
real-world instances of the Docklands Light Railway with 108--336 active
trains, the model predicts reductions of 19.3\%--28.8\%, solving each
instance in under 0.11 seconds. Our model is being integrated into Hitachi
Rail Canada's industrial timetable compiler.

\end{abstract}


\section{Introduction}
\paragraph{Background.}

The United Nations projects that the share of the world's population living in urban areas, based on national definitions, will rise from approximately 58\% in 2025 to 67\% in 2050 \citep{unitednations2025}, placing sustained pressure on existing transportation infrastructure. Rail transit, with its high passenger-carrying capacity and strong safety record \citep{Abril2008}, is well positioned to absorb much of this growth. Communications-based train control (CBTC) systems play a central role in the operation of modern metro networks \citep{pascoe2009communication}. Used by some of the world's busiest metros, they provide high-resolution train location determination and continuous data communication, which improve both safety and efficiency and allow operational schedules to be adjusted dynamically in response to congestion and fluctuating demand. Because CBTC systems can implement automatic driving strategies and execute a prescribed timetable precisely \citep{IRSE2023}, they are particularly well suited to energy-efficient operation. These capabilities make fast and energy-efficient timetable optimization a natural way to improve both the efficiency and the sustainability of CBTC-enabled metro networks.

\begin{figure}
\begin{centering}
\includegraphics[scale=0.4]{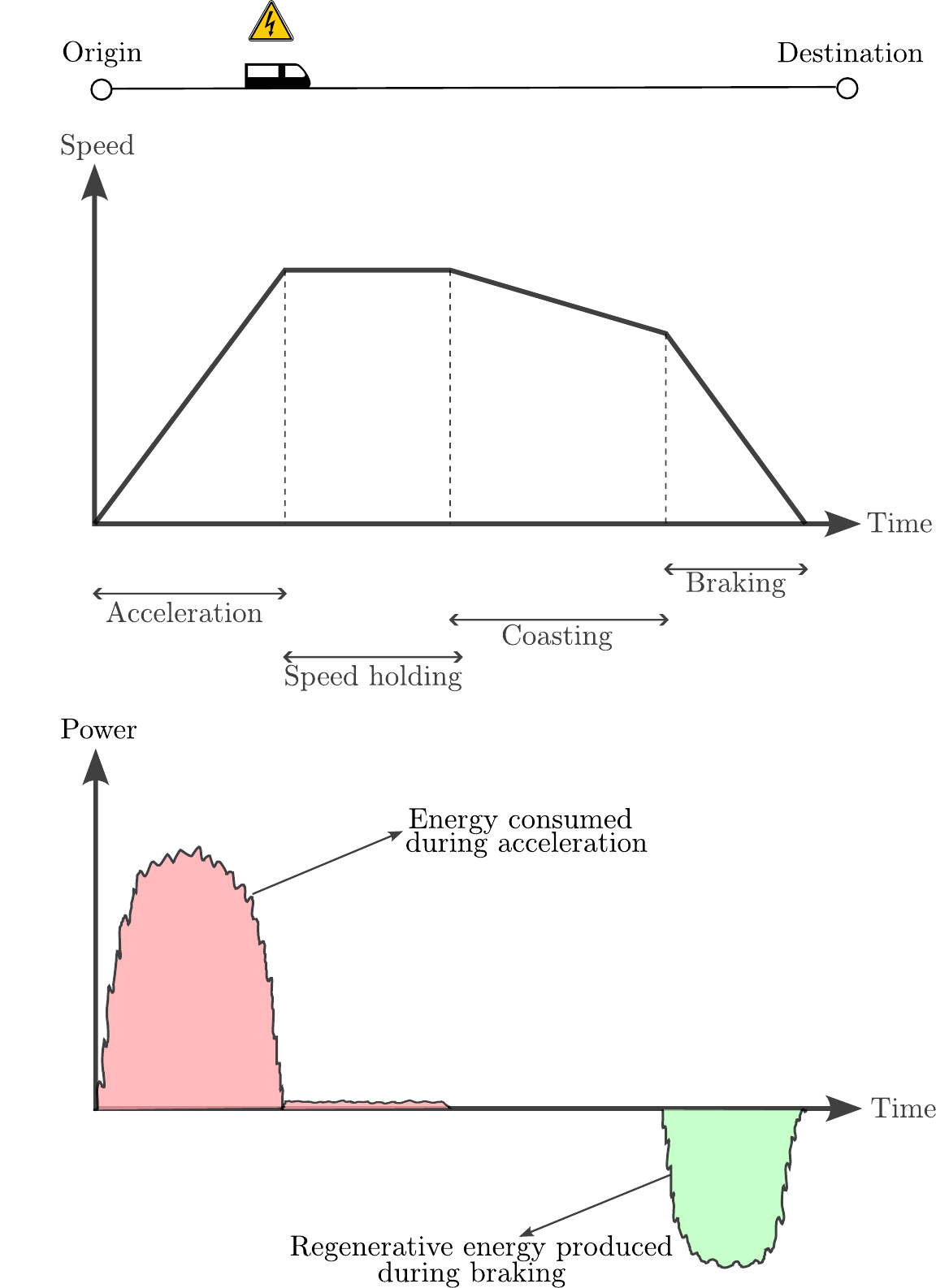}
\par\end{centering}
\caption{Illustration of a train's energy-use cycle under CBTC operation. For a trip between two platforms, the qualitative speed profile consists of acceleration, speed holding, coasting, and braking. The lower plot shows the corresponding power profile: the acceleration interval (red) accounts for most traction energy demand, whereas the braking interval (green) produces regenerative energy. With suitable scheduling, that recovered energy can be supplied to nearby accelerating trains.\label{Fig:SpeedProfileOfTrain}}
\end{figure}

\paragraph{Motivation.}

In CBTC-enabled metro systems, trains are electrically powered. During a trip from one platform to the next, the train follows the four-phase speed pattern shown in Figure \ref{Fig:SpeedProfileOfTrain}: acceleration, speed holding, coasting, and braking \citep{Howlett1995}. Most of the energy consumptions occurs in acceleration phase. During braking, the traction motors operate as generators and produce \textit{regenerative braking energy}. Because CBTC systems can implement a given railway timetable with high degree of accuracy to the trains operating in the railway network, this energy can be transferred from braking trains to nearby trains that are accelerating. For perspective, the New York City transit railway system consumes more than 1,600 GWh of electricity per year, and all trains installed since 2018 are capable of producing and transferring regenerative braking energy; under favorable conditions, the recovered energy may amount to as much as 50\% of the consumed energy \citep{mohamed2018white}.

\paragraph{Energy-optimal timetable.}

An energy-optimal timetable offers substantial energy savings without requiring infrastructural alterations to a CBTC-enabled railway system. Formally, a railway timetable is a data structure containing the arrival and departure times of each train at every platform visited during a designated service period, typically spanning 18 to 24 hours in most operational networks. An energy-optimal timetable chooses these temporal decision variables to minimize the network's \textit{effective energy consumption}: the \textit{total energy consumption of all trains during their accelerating phases} minus \textit{the total regenerative braking energy successfully transferred to the accelerating trains from eligible braking counterparts}.

\paragraph{Related work.}

Energy-efficient railway operation is studied at two levels: energy-efficient train control, which shapes a single train's speed profile between two stops, and energy-efficient timetabling, which coordinates the event times of many trains. \citet{Scheepmaker2017} review both levels and their interaction, \citet{Yang2016survey} survey the urban-rail line of work, and \citet{GonzalezGil2013} review the technologies for recovering regenerative braking energy. One of the technologies is the train-to-train transfer exploited in this paper, which requires no dedicated storage hardware.

Within the timetabling level, work on regenerative-energy utilization spans several modeling paradigms. \citet{Pena-Alcaraz2012} study Madrid Metro Line 3 using a mixed-integer programming formulation tailored to a single line. \citet{Yang2013} initiated a cooperative-scheduling line that maximizes the overlap time between the braking and accelerating phases of successive trains (the objective that our energy-weighted overlaps refine), and \citet{wang2022improving} design a mixed-integer linear programming model to increase the synchronization between accelerating and braking trains at the same station; \citet{Su2013} couple the timetable with energy-efficient driving regimes. \citet{Li2014} use a genetic algorithm to search for timetables that balance regenerative-energy utilization with train tractive energy; \citet{Le2014} and \citet{liu2018energy} develop a nonlinear integer program solved by simulated annealing and by a hybrid algorithm, respectively. \citet{DasGuptaACC15} give a more computationally tractable mixed-integer program for coordinating regenerative-energy exchange among suitable train pairs and demonstrate it on the Docklands Light Railway, but their model does not optimize actual energy savings directly. \citet{Yin2017} embed regenerative-energy utilization and dynamic passenger demand in mixed-integer linear programs solved by Lagrangian relaxation, and \citet{Canca2017} design energy-efficient timetables for two-way rapid transit lines. The bi-objective model of \citet{yang2020bi} combines timetable optimization with energy allocation and passenger assignment through a non-dominated sorting genetic algorithm \citep[see also][]{yang2021efficient}; \citet{Gemander2023} optimize underground timetables under uncertain dwell and running times through stochastic programming; and \citet{barmann2021efficient} design exact Benders decomposition schemes to optimize train timetables for peak power reduction and regenerative-energy use improvement.
A common pattern in these mixed-integer formulations above is a set of binary
variables encoding whether the braking phase of one train overlaps the accelerating
phase of another; in Section~\ref{sec:regen_descrip} we show that this overlap can be modeled
as a concave function of the event times, so no integer variables are actually needed.

A complementary line of work combines timetabling with train control. The model in \citet{aLi2014a} jointly optimize dynamic train scheduling and speed control, \citet{su2019energy} and \citet{zhang2024energy} both use a two-level optimization approach: speed profile is optimized at the first level, while the optimal timetable is solved given the optimal speed profile at the second level.
\citet{Luan2018} extend integrated real-time traffic management to energy-efficient train operations, and \citet{wang2017multi} address multi-train trajectory optimization through a multiple-phase optimal-control formulation solved by a pseudospectral method, later adapting it to retime a fixed timetable for improved energy efficiency \citep{wang2019multi}. \citet{shi2026integrated} develop a decomposition approach to solve the integrated mixed-integer linear programming model for train rescheduling and speed control during disruption, and \citet{howlett2023optimal} formulate an unconstrained convex model to solve for optimal timetable and driving strategies for each train in a fleet. In contrast, we focus on building our optimization model based on how CBTC system works in practice: the on-board controller computes the energy-optimal speed profile for whatever trip times are assigned, and the optimization chooses only the event times, which is what keeps the network-wide model a single linear program.

From the standpoint of deployed industrial tools, the two-stage framework of \citet{gupta2016two} appears to be, to our knowledge, the only earlier model embedded in a CBTC timetable compiler. That work computes an energy-efficient timetable through two consecutive linear programs: the first stage minimizes network-wide energy use over all trains subject to the railway constraints, and the second stage fixes the resulting trip times and then maximizes regenerative-energy transfer between suitable train pairs.  However, this model in \citet{gupta2016two}  has some limitations that limit its real-world applicability. First, although each stage's linear program solves quickly, the handoff between the stages requires an intermediate simulation step that can take hours, so the total runtime is much longer. Second, the two-stage objective does not model actual network energy consumption directly; instead, it optimizes proxy quantities whose practical effect on energy usage must be assessed afterward by extensive physics-based simulation, which limits managerial usefulness when explicit energy predictions are desired. Third, the second stage relies on the exact locations of power-consumption peaks, information that is often unavailable in practice. These limitations in \citet{gupta2016two} have motivated us to design the new optimization model for CBTC systems that we describe in this paper.

Relative to that model, the present one keeps the operational constraint families of Section~\ref{sec:Modelling-the-constraints} and the dwell-midpoint pairing of braking with accelerating trains, and redesigns the optimization itself: transferred regenerative energy is now fitted directly in energy units as an affine function of phase overlap, the effective phases vary affinely with trip time, the overlaps are captured exactly through linear hypograph constraints in place of the earlier $\ell_{1}$ alignment heuristic, and a single linear program minimizes (and predicts, without physics-based simulation) the network's effective energy consumption.

\paragraph{Contribution.}

We develop a single-stage linear optimization model that constructs energy-optimal timetables for all trains in a CBTC-enabled metro railway network over an entire service period. Our contributions are as follows.
\begin{itemize}
\item \emph{A single-stage linear model.} The model minimizes effective energy consumption (traction energy minus transferred regenerative energy) in one linear program, in contrast with earlier approaches that rely on NP-hard mixed-integer formulations or on multiple stages interleaved with time-consuming simulation.
\item \emph{Data-driven, simulation-free energy prediction.} Traction energy, transferred regenerative energy, and their interaction are modeled from operational or simulated data, so the model predicts the effective energy consumption of a timetable without invoking a physics-based simulator. This makes the model suitable for managerial settings in which explicit energy predictions are needed quickly.
\item \emph{A convex-analytic model of regenerative-energy transfer.} We show that the
 signed overlap between the effective braking and accelerating phases of a train pair
 is concave in the timetable, being a pointwise minimum of affine functions of the
 event times. Maximizing the transferred regenerative energy is therefore a convex
 problem, which a hypograph reformulation expresses through linear constraints without
 loss of optimality. 
Prior work either captures this synchronization structure with binary variables or other nonconvex constraints, or relies on multiple stages interleaved with expensive simulation.
\item \emph{Validation on two real-world networks and industrial integration.} On 11 full-day operational instances of Shanghai Metro Line 8 (one of the busiest metro services in the world) with 1,000--1,332 active trains, the model solves in under 0.7 seconds per instance on a laptop and predicts effective-energy reductions of 19.52\%--28.17\% relative to the existing operational timetables; these predictions are corroborated by Hitachi Rail's physics-based simulator \mbox{\textsf{SPSIM}}. On 17 instances of the Docklands Light Railway benchmark with 108--336 active trains, the model predicts reductions of 19.33\%--28.82\%. The model is being integrated into Hitachi Rail Canada's industrial timetable compiler.
\end{itemize}

This work arose from a collaboration with Hitachi Rail Canada, the largest provider of CBTC systems worldwide. The collaboration began in 2023 with Thales Ground Transportation Systems (GTS). On May~31, 2024, Hitachi Rail completed its acquisition of Thales's worldwide GTS business for \texteuro 1.66~billion \citep{hitachiGTS2024}. Due to this transaction, the established Canadian GTS organization became Hitachi Rail Canada, and our collaboration continued with the Toronto-based team under Hitachi Rail Canada. We therefore use the current name \emph{Hitachi Rail Canada} throughout this paper.

\subsection{Organization}

The remainder of the paper is structured as follows. Section~\ref{notationAndNotions} introduces notation and preliminary concepts. Section~\ref{sec:Modelling-the-constraints} presents the operational constraints required for a metro railway network. In Section~\ref{sec:Modeling-the-objective}, we develop the effective-energy-consumption objective. Section~\ref{sec:Final-optimization-model} assembles the full optimization model. Section~\ref{Numerical_Study} reports experiments on Shanghai Metro Line~8 and the Docklands Light Railway benchmark, and Section~\ref{sec:archframe} describes the architecture used for industrial deployment of the model.

\section{Notation and preliminaries \label{notationAndNotions}}

Unless stated otherwise, every set appearing in this paper is finite and equipped with a strict order. For such a set $S$, we write $|S|$ for its cardinality and $S(i)$ for its $i$-th element.

\begin{longtable}[c]{ll}
\caption{List of symbols in the order they appear in the paper.}
\tabularnewline
{\scriptsize{}\textbf{Symbol} } & {\scriptsize{}{}\textbf{Description}}
\tabularnewline
\endfirsthead
\multicolumn{2}{l}{{\scriptsize{}\emph{Table~\thetable{} (continued)}}}
\tabularnewline
{\scriptsize{}\textbf{Symbol} } & {\scriptsize{}{}\textbf{Description}}
\tabularnewline
\endhead
{\scriptsize{}$\mathcal{N}$ } & {\scriptsize{}{}All platforms in the railway network }\tabularnewline
{\scriptsize{}{}$\mathcal{A}$ } & {\scriptsize{}{}All tracks in the railway network}\tabularnewline
{\scriptsize{}{}$\mathcal{T}$ } & {\scriptsize{}{}Collection of trains under consideration }\tabularnewline
{\scriptsize{}{}$\mathcal{N}^{t}$ } & {\scriptsize{}{}Platforms visited by train $t$, listed in chronological order }\tabularnewline
{\scriptsize{}{}$\mathcal{A}^{t}$ } & {\scriptsize{}{}Tracks traversed by train $t$, listed in chronological order }\tabularnewline
{\scriptsize{}{}$a_{i}^{t}$ } & {\scriptsize{}{}Decision variable for the arrival time of train $t$ at platform $i$ }\tabularnewline
{\scriptsize{}{}$d_{i}^{t}$ } & {\scriptsize{}{}Decision variable for the departure time of train $t$ from platform $i$ }\tabularnewline
{\scriptsize{}{}$[\underline{\tau}_{ij}^{t},\overline{\tau}_{ij}^{t}]$ } & {\scriptsize{}{}Admissible trip-time window for train $t$ from platform $i$ to platform $j$ on a track }\tabularnewline
{\scriptsize{}$\varphi$} & {\scriptsize{}Set of crossing-overs in the network }\tabularnewline
{\scriptsize{}{}$\mathcal{B}_{ij}$ } & {\scriptsize{}{}Train pairs associated with turn-around movements on crossing-over $(i,j)$ }\tabularnewline
{\scriptsize{}{}$[\underline{\kappa}_{ij}^{tt'},\overline{\kappa}_{ij}^{tt'}]$ } & {\scriptsize{}{}Admissible trip-time window for the train pair $(t,t')$ on crossing-over $(i,j)$ }\tabularnewline
{\scriptsize{}{}$[\underline{\delta}_{i}^{t},\overline{\delta}_{i}^{t}]$ } & {\scriptsize{}{}Allowable dwell-time window for train $t$ at platform $i$ }\tabularnewline
{\scriptsize{}{}$\chi$ } & {\scriptsize{}{}Platform pairs located at the same interchange station }\tabularnewline
{\scriptsize{}{}$\mathcal{C}_{ij}$ } & {\scriptsize{}{}Connecting train pairs associated with platform pair $(i,j)\in\chi$ }\tabularnewline
{\scriptsize{}{}$[\underline{\chi}_{ij}^{tt'},\overline{\chi}_{ij}^{tt'}]$ } & {\scriptsize{}{}Transfer-time window between train $t$ at platform $i$ and train $t'$ at platform $j$ }\tabularnewline
{\scriptsize{}{}$\mathcal{H}_{ij}$ } & {\scriptsize{}{}Ordered train pairs that successively traverse track $(i,j)$ }\tabularnewline
{\scriptsize{}{}$h_{i}^{tt'}$ } & {\scriptsize{}{}Required headway between the departures of trains $t$ and $t'$ from platform $i$ }\tabularnewline
{\scriptsize{}{}$[\underline{\tau}_{\mathcal{P}}^{t},\overline{\tau}_{\mathcal{P}}^{t}]$ } & {\scriptsize{}{}Allowable total travel-time window for train $t$ along its train path} \tabularnewline
{\scriptsize{}$E_{i,j,t}^{\textrm{con,tr}}$} & {\scriptsize{}Traction energy consumed by train $t$ while accelerating on the trip from platform $i$ to platform $j$}\tabularnewline
{\scriptsize{}$E_{i,j,t,t^{\prime}}^{\textrm{con,cr}}$} & {\scriptsize{}Traction energy consumed by train $t$ while accelerating through crossing-over $(i,j)\in\varphi$}\tabularnewline
{\scriptsize{}{}$\Omega$ } & {\scriptsize{}{}Feasible platform pairs for regenerative-energy transfer}\tabularnewline
{\scriptsize{}{}$\mathcal{T}_{i}$ } & {\scriptsize{}{}Trains that arrive at, dwell at, and depart from platform $i$}\tabularnewline
{\scriptsize{}$\mu_{i}^{t}$ } & {\scriptsize{}{}Dwell midpoint $(\bar{a}_{i}^{t}+\bar{d}_{i}^{t})/2$ of train $t$ at platform $i$ in the reference timetable}\tabularnewline
{\scriptsize{}{}$\overset{\rightharpoonup}{t}$ } & {\scriptsize{}{}Train temporally close to the right of $t$}\tabularnewline
{\scriptsize{}{}$\overset{\leftharpoonup}{t}$ } & {\scriptsize{}{}Train temporally close to the left of $t$}\tabularnewline
{\scriptsize{}{}$\tilde{t}$ } & {\scriptsize{}{}Train temporally close to $t$}\tabularnewline
{\scriptsize{}{}$\mathcal{E}$ } & {\scriptsize{}{}Synchronization events between suitable train pairs}\tabularnewline
{\scriptsize{}{}$\overset{\rightharpoonup}{\mathcal{E}}$ } & {\scriptsize{}{}Subset of $\mathcal{E}$ containing tuples of the form $(i,j,t,\overset{\rightharpoonup}{t})$}\tabularnewline
{\scriptsize{}{}$\overset{\leftharpoonup}{\mathcal{E}}$ } & {\scriptsize{}{}Subset of $\mathcal{E}$ containing tuples of the form $(i,j,t,\overset{\leftharpoonup}{t})$}\tabularnewline
{\scriptsize{}$E_{i,j,t,\overset{\rightharpoonup}{t}}^{\textrm{reg}}$} & {\scriptsize{}Regenerative energy transferred to train $t$ (accelerating at $i$) by train $\overset{\rightharpoonup}{t}$ (braking at $j$)}\tabularnewline
{\scriptsize{}$\sigma_{ij}^{t\overset{\rightharpoonup}{t}}$} & {\scriptsize{}Overlap of the accelerating phase of $t$ at $i$ with the braking phase of $\overset{\rightharpoonup}{t}$ at $j$}\tabularnewline
{\scriptsize{}$E_{i,j,t,\overset{\leftharpoonup}{t}}^{\textrm{reg}}$} & {\scriptsize{}Regenerative energy transferred to train $\overset{\leftharpoonup}{t}$ (accelerating at $j$) by train $t$ (braking at $i$)}\tabularnewline
{\scriptsize{}$\sigma_{ij}^{t\overset{\leftharpoonup}{t}}$} & {\scriptsize{}Overlap of the accelerating phase of $\overset{\leftharpoonup}{t}$ at $j$ with the braking phase of $t$ at $i$}\tabularnewline
{\scriptsize{}$\underline{\beta}_{i}^{t},\overline{\beta}_{i}^{t}$ } & {\scriptsize{}{}Start and end offsets of the effective braking phase of train $t$ before its arrival at platform $i$ ($\underline{\beta}_{i}^{t}\geq\overline{\beta}_{i}^{t}$)}\tabularnewline
{\scriptsize{}$\underline{\alpha}_{i}^{t},\overline{\alpha}_{i}^{t}$ } & {\scriptsize{}{}Start and end offsets of the effective accelerating phase of train $t$ after its departure from platform $i$}\tabularnewline
{\scriptsize{}$\varpi_{ij}^{t\tilde{t}},\psi_{ij}^{t\tilde{t}}$ } & {\scriptsize{}{}Auxiliary variables that linearize the overlap $\sigma_{ij}^{t\tilde{t}}$ in the final model \eqref{eq:final_lp_model}}\tabularnewline
\end{longtable}

Platforms are elements of the set $\mathcal{N}$. We say that two platforms are \textit{opposite} if they belong to the same station but trains on those platforms travel in opposite directions. Two tracks are likewise opposite if they join the same pair of stations in opposite directions. We use the term \textit{track} for a directed arc joining two distinct, non-opposite platforms, and collect all such arcs in $\mathcal{A}$.

The railway network is represented by the directed graph $\mathcal{G}=(\mathcal{N},\mathcal{A})$. A \textit{train-line} is a directed path whose nodes correspond to non-opposite platforms and whose arcs correspond to non-opposite tracks. Next, a \textit{crossing-over} is a directed arc linking two train-lines at their terminal platforms. In CBTC terminology, when a physical train reaches the terminal platform of one line, traverses a crossing-over, and then begins service on another line, the CBTC management treats the movement \emph{operationally} involving two distinct trains (even though physically there is one physical train). The set of trains considered in the model is $\mathcal{T}$. For any $t\in\mathcal{T}$, the platforms and tracks visited by $t$, listed in temporal order, are denoted by $\mathcal{N}^{t}\subseteq\mathcal{N}$ and $\mathcal{A}^{t}\subseteq\mathcal{A}$, respectively. The directed path formed by these ordered platforms and tracks is called the \textit{train-path} of $t$.

The optimization variables are the arrival and departure times of trains at the platforms they visit. We write $a_{i}^{t}$ for the arrival time of train $t\in\mathcal{T}$ at platform $i\in\mathcal{N}^{t}$, and $d_{j}^{t}$ for the departure time of train $t$ from platform $j\in\mathcal{N}^{t}$. 

\section{Modeling the operational constraints \label{sec:Modelling-the-constraints}}

The operational constraints relate the arrival and departure events of the trains and are necessary for the proper functioning of the railway network; we present them in Sections~\ref{subsec:Trip-time-constraint}--\ref{subsec:Total-travel-time}. First, we describe how the model protects against uncertainty in the constraint data.

\paragraph{Incorporating robust constraints through box uncertainty.}\label{subsec:Robust-constraints-through}

\begin{figure}
\begin{centering}
\includegraphics[scale=1.75]{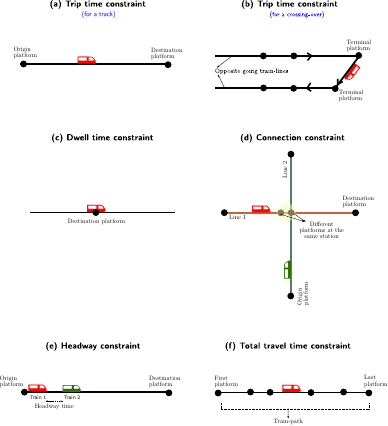}
\par\end{centering}
\caption{Graphical summary of the operational constraints that any feasible timetable in a metro railway network must satisfy. \label{Fig:AllCon}}
\end{figure}

In practice, there can be some uncertainty in the data associated with the operational constraints, which can cause them to diverge from their intended values. To guard against such uncertainties, we integrate box uncertainty into our model. Abstractly, every operational constraint below can be represented as $\ell\leq x-y\leq u$, where $x$ and $y$ stand for decision variables, and $\ell$ and $u$ refer to the problem data. We model the uncertainty in the data by the boxes $\ell\in[\ell_{\textrm{lb}},\ell_{\textrm{ub}}]$ and $u\in[u_{\textrm{lb}},u_{\textrm{ub}}]$. Requiring $\ell\leq x-y\leq u$ to hold for every $\ell$ and $u$ in these boxes amounts to imposing the worst case in each box, so the robust counterpart of the constraint is $\ell_{\textrm{ub}}\leq x-y\leq u_{\textrm{lb}}$. For all our instances, $\ell_{\textrm{ub}}\leq u_{\textrm{lb}}$ holds; therefore, each robust window is nonempty. For brevity, we present the robust versions of our constraints directly in the representations that follow. While all the constraints share this abstract form, they represent very different physical requirements, which Figure~\ref{Fig:AllCon} illustrates.

\subsection{Trip-time constraint \label{subsec:Trip-time-constraint}}

Trip-time limits are central because they govern both traction energy use and also regenerative-energy production. As we discussed earlier, nearly all traction energy for a train's trip from one platform to the next is spent during acceleration, while regenerative energy is generated during braking. We therefore distinguish the following two trip-time constraints.

\subsubsection{Track-based trip-time constraint}

Consider train $t\in\mathcal{T}$ traveling from platform $i$ to platform $j$ on track $(i,j)\in\mathcal{A}^{t}$. Let $d_{i}^{t}$ and $a_{j}^{t}$ denote its departure and arrival times. The resulting travel time must lie in the admissible window $[\underline{\tau}_{ij}^{t},\overline{\tau}_{ij}^{t}]$. Therefore,
\begin{align}
\underline{\tau}_{ij}^{t}\leq a_{j}^{t}-d_{i}^{t}\leq\overline{\tau}_{ij}^{t}, \quad \textrm{for } (i,j)\in\mathcal{A}^{t}, \, t\in\mathcal{T}.  & \tag{\texttt{TRACK}}\label{eq:TripTimeConstraint}
\end{align}
This constraint is shown in Figure~\ref{Fig:AllCon}(a).

\subsubsection{Crossing-over-based trip-time constraint}

A crossing-over is a directed arc connecting two train-lines at their terminal platforms. Consider $(i,j)\in\varphi$, where platforms $i$ and $j$ belong to different lines, and let $\mathcal{B}_{ij}$ denote the train pairs associated with turn-around movements on this crossing-over. For $(t,t')\in\mathcal{B}_{ij}$, train $t$ departs platform $i$, traverses $(i,j)$, and then continues from platform $j$ under the new operational label $t'$. The elapsed time between departure from $i$ as train $t$ and arrival at $j$ as train $t'$ must remain in the admissible window $[\underline{\kappa}_{ij}^{tt'},\overline{\kappa}_{ij}^{tt'}]$. Hence,
\begin{align}
\underline{\kappa}_{ij}^{tt'}\leq a_{j}^{t'}-d_{i}^{t}\leq\overline{\kappa}_{ij}^{tt'}, \quad \textrm{for } (t,t')\in\mathcal{B}_{ij}, \, (i,j)\in\varphi. & \tag{\texttt{CROSS}}\label{eq:turnAroundConstraint}
\end{align}

This constraint is shown in Figure~\ref{Fig:AllCon}(b).

\subsection{Dwell-time constraint\label{subsec:Dwell-time-constraint}}
When train $t\in\mathcal{T}$ reaches platform $i\in\mathcal{N}^{t}$, it remains there for a dwell interval in $[\underline{\delta}_{i}^{t},\overline{\delta}_{i}^{t}]$ so that old passengers can get off and new passengers can board. Accordingly, the difference between departure and arrival times at that platform must satisfy
\begin{align}
\underline{\delta}_{i}^{t}\leq d_{i}^{t}-a_{i}^{t}\leq\overline{\delta}_{i}^{t}, \quad \textrm{for } i\in\mathcal{N}^{t}, \, t\in\mathcal{T}. & \tag{\texttt{DWELL}}\label{eq:DwellTimeConstraint}
\end{align}
Figure~\ref{Fig:AllCon}(c) depicts this constraint. The index $i$ ranges over all of $\mathcal{N}^{t}$ because every train dwells at each platform on its path, including the first and last ones: the train reaches $\mathcal{N}^{t}(1)$ either from the depot or via a turn-around movement, and after departing $\mathcal{N}^{t}(|\mathcal{N}^{t}|)$ it either returns to the depot or begins service on another line through a turn-around. 

\subsection{Connection constraint \label{subsec:Connection-constraint}}

Passengers may need to transfer at interchange stations when no direct train serves their desired origin--destination pair. Let $\chi\subseteq\mathcal{N}\times\mathcal{N}$ collect the platform pairs at which such transfers are possible. For $(i,j)\in\chi$, the two platforms belong to the same station. Let $\mathcal{C}_{ij}$ denote the train pairs $(t,t')$ for which passengers arrive on train $t$ at platform $i$ and continue on train $t'$ departing from platform $j$. In a feasible timetable, the time available for this transfer has to remain inside the prescribed window. Accordingly, the connection constraint is
\begin{align}
\underline{\chi}_{ij}^{tt'}\leq d_{j}^{t'}-a_{i}^{t}\leq\overline{\chi}_{ij}^{tt'}, \quad \textrm{for } (t,t')\in\mathcal{C}_{ij}, \, (i,j)\in\chi, \tag{\texttt{CONNECT}}\label{eq:ConnectionConstraint}
\end{align}
where $\underline{\chi}_{ij}^{tt'}$ and $\overline{\chi}_{ij}^{tt'}$ are the lower and upper limits of the transfer window. Figure~\ref{Fig:AllCon}(d) illustrates this constraint.

\subsection{Headway constraint \label{subsec:Headway-constraint}}

Railway operations also require a minimum temporal separation between consecutive trains, called the \textit{headway}. For track $(i,j)\in\mathcal{A}$, let $\mathcal{H}_{ij}$ denote the ordered train pairs that traverse the track successively. If $(t,t')\in\mathcal{H}_{ij}$, then train $t'$ follows train $t$ on the same track. The headway requirements at platforms $i$ and $j$ are $h_{i}^{tt'}$ and $h_{j}^{tt'}$, which leads to the following constraint: 
\begin{align}
h_{i}^{tt'}\leq d_{i}^{t'}-d_{i}^{t}\;\textrm{ and }\ h_{j}^{tt'}\leq d_{j}^{t'}-d_{j}^{t}, \quad \textrm{for } (t,t')\in\mathcal{H}_{ij}, \, (i,j)\in\mathcal{A}. & \tag{\texttt{HEADWAY}}\label{eq:SafetyConstraint1}
\end{align}
The headway constraint is shown in Figure~\ref{Fig:AllCon}(e). 

The headway constraints also restrict regenerative-energy synchronization. When a train accelerates away from a platform, that platform cannot be occupied immediately by another train, so pairing two consecutive trains on the same track (one accelerating and one braking at the same platform) is never feasible from a safety standpoint. We account for this restriction when constructing suitable train pairs in Section~\ref{sec:regen_descrip}.

\subsection{Total travel-time constraint \label{subsec:Total-travel-time}}

Reliable service requires the total travel time of each train $t\in\mathcal{T}$
to lie within a specified window $[\underline{\tau}_{\mathcal{P}}^{t},\overline{\tau}_{\mathcal{P}}^{t}]$, where the subscript $\mathcal{P}$ abbreviates \emph{train-path}:
the time elapsed between the train's departure from its first
platform $\mathcal{N}^{t}(1)$ and its arrival at its last platform
$\mathcal{N}^{t}(|\mathcal{N}^{t}|)$ must satisfy
\begin{align}
\underline{\tau}_{\mathcal{P}}^{t}\leq a_{\mathcal{N}^{t}(|\mathcal{N}^{t}|)}^{t}-d_{\mathcal{N}^{t}(1)}^{t}\leq\overline{\tau}_{\mathcal{P}}^{t}, \quad \textrm{for } t\in\mathcal{T}. & \tag{\texttt{TRAVEL}}\label{eq:TotalTravelTimeConstraints}
\end{align}
Keeping total travel times within these windows protects connections and preserves network capacity. This constraint is shown in Figure~\ref{Fig:AllCon}(f).

\subsection*{Domain of the event times}

Bounds on event times are set relative to the service horizon. We normalize the beginning of the service period to time zero, which corresponds to the day's first departure from some platform. An upper bound on the final event can then be obtained by assigning every trip time and dwell time its maximal allowed value; call this bound $m>0$. The arrival and departure variables are therefore restricted by
\begin{align}
0\leq a_{i}^{t}\leq m,\;0\leq d_{i}^{t}\leq m, \quad \textrm{for } i\in\mathcal{N}^{t},\, t\in\mathcal{T}. & \tag{\texttt{DOMAIN}}\label{eq:domain}
\end{align}

\section{Modeling the effective energy consumption \label{sec:Modeling-the-objective}}

This section develops the objective representing \emph{effective energy consumption}. We define this quantity as the total traction energy consumed during train acceleration, $E^{\textrm{con}}$, minus the regenerative braking energy transferred to eligible accelerating trains, $E^{\textrm{reg}}$. Hence the objective is $E^{\textrm{con}}-E^{\textrm{reg}}$. The data-driven model for $E^{\textrm{con}}$ is presented in Section~\ref{sec:consum}, followed by the model for $E^{\textrm{reg}}$ in Section~\ref{sec:regen_descrip}.

\subsection{Modeling total consumed energy $E^{\textrm{con}}$}\label{sec:consum}

\paragraph{Relationship between energy consumption and the trip-time constraint.}

Most electrical energy used by a train is drawn during acceleration between successive platforms, so trip-time decisions are a primary determinant of both traction consumption and regenerative-energy generation. Once a trip time is specified, CBTC tools can compute the corresponding energy-optimal speed profile. In our study we use \textsf{KDC} \citep{selTrac}, a proprietary physics-based train-performance simulator used within Hitachi Rail's CBTC toolchain, and it computes speed profiles that follow the acceleration, speed-holding, coasting, and braking sequence. This profile is consistent with the theoretical optimum described by \citet{Howlett1995}; see \citet{Jiaxin1993,Howlett2000,Howlett2009,Khmelnitsky2000,Liu2003} for foundational treatments and \citet{Albrecht2016} for a broader review. The speed profile determines the train's electrical power draw and regeneration over time, and therefore its associated power-versus-time curve, illustrated qualitatively in Figure~\ref{Fig:SpeedProfileOfTrain}.

\paragraph{Energy consumption of all trains.} We represent the network-wide traction energy used by all trains as 

\[
E^{\textrm{con}}\triangleq\sum_{t\in\mathcal{T}}\sum_{(i,j)\in\mathcal{A}^{t}}\underbrace{E_{i,j,t}^{\textrm{con,tr}}}_{\substack{\textrm{depends on }\\
\textrm{trip time }(a_{j}^{t}-d_{i}^{t})
}
}\quad+\sum_{(i,j)\in\varphi,(t,t^{\prime})\in\mathcal{B}_{ij}}\underbrace{E_{i,j,t,t^{\prime}}^{\textrm{con,cr}}}_{\substack{\textrm{depends on }\\
\textrm{trip time }(a_{j}^{t'}-d_{i}^{t})
}
},
\]
where $E_{i,j,t}^{\textrm{con,tr}}:\mathbb{R}_{++}\to\mathbb{R}_{++}$ denotes the traction energy used by train $t\in\mathcal{T}$ while accelerating on the trip from platform $i$ to platform $j$ for $(i,j)\in\mathcal{A}^{t}$; in a CBTC setting its argument is the trip time $a_{j}^{t}-d_{i}^{t}$. Likewise, $E_{i,j,t,t^{\prime}}^{\textrm{con,cr}}:\mathbb{R}_{++}\to\mathbb{R}_{++}$ represents the traction energy used by train $t$ when it traverses crossing-over $(i,j)\in\varphi$ and then continues from platform $j$ under label $t^{\prime}$, with $(t,t^{\prime})\in\mathcal{B}_{ij}$. Its argument is the corresponding crossing-over travel time $a_{j}^{t'}-d_{i}^{t}$.

\begin{figure}
\begin{centering}
\includegraphics[scale=0.8]{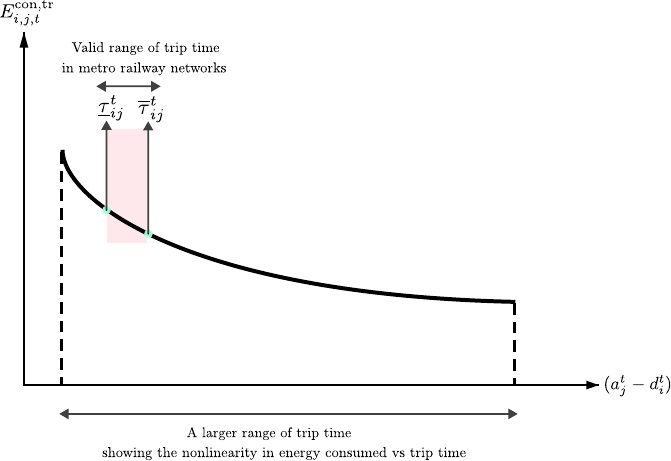}
\par\end{centering}
\caption{Qualitative behavior of the traction-energy function for train $t$ on track $(i,j)$ as the trip time $a_j^t - d_i^t$ varies. The wide horizontal range is included only to display the nonlinear shape observed in long inter-city runs. In metro operation, the feasible interval $[\underline{\tau}_{ij}^{t}, \overline{\tau}_{ij}^{t}]$ is typically measured in seconds, and over that narrow range an affine approximation of consumed energy is reasonable. \label{Fig:ConsumedEngTrain}}
\end{figure}

\paragraph{Data-driven approach to model the energy consumption.}

The functions $E_{i,j,t}^{\textrm{con,tr}}$ and $E_{i,j,t,t^{\prime}}^{\textrm{con,cr}}$ generally do not admit a tractable closed form \citep{Howlett2009, gupta2016two}, although their qualitative behavior is well understood. Under an optimal speed profile, each is monotone nonincreasing in trip time \citep{Milroy1980}. Empirical studies indicate that the same monotonic trend persists on average even under manual or otherwise suboptimal driving; see \citet[ Figure~\ref{Fig:SpeedProfileOfTrain}]{pena2011approximate}. Moreover, the relevant energy data are accessible in practice, either from historical operations or from physics-based simulation. A key feature of CBTC metro systems is that the permissible variation in the trip-time constraints \eqref{eq:TripTimeConstraint} and \eqref{eq:turnAroundConstraint} is usually only a few seconds. We therefore adopt the following assumption:

\begin{assumption} \label{assum:triptime} For every track-based or crossing-over trip-time constraint, the width of the admissible time window is on the order of seconds; equivalently, $\overline{\tau}_{ij}^{t}-\underline{\tau}_{ij}^{t}$ in \eqref{eq:TripTimeConstraint} and $\overline{\kappa}_{ij}^{tt'} - \underline{\kappa}_{ij}^{tt'}$ in \eqref{eq:turnAroundConstraint} are on the order of seconds. \end{assumption}

Assumption \ref{assum:triptime}, together with the monotone dependence on trip time, justifies affine approximations of the energy-consumption functions. Figure~\ref{Fig:ConsumedEngTrain} illustrates this idea for a track-based trip; the crossing-over case behaves analogously. To obtain the best affine fit, we first apply least-squares to energy-versus-trip-time data, using either historical observations or outputs from physics-based simulators such as the \textsf{KDC} simulator employed in this work \citep{selTrac}.

In other words, we approximate consumed energy as
\begin{equation}
\begin{alignedat}{1}E_{i,j,t}^{\textrm{con,tr}}&\;\triangleq\; c_{i,j,t}^{\textrm{con,tr}}(a_{j}^{t}-d_{i}^{t})+b_{i,j,t}^{\textrm{con,tr}},\quad \textrm{for }(i,j)\in\mathcal{A}^{t},\,t\in\mathcal{T},\\
E_{i,j,t,t^{\prime}}^{\textrm{con,cr}} & \;\triangleq\; c_{i,j,t,t^{\prime}}^{\textrm{con,cr}}(a_{j}^{t'}-d_{i}^{t})+b_{i,j,t,t^{\prime}}^{\textrm{con,cr}},\quad \textrm{for }(i,j)\in\varphi,(t,t^{\prime})\in\mathcal{B}_{ij},
\end{alignedat}
\label{eq:eng_con_model}
\end{equation}
where $(c_{i,j,t}^{\textrm{con,tr}},b_{i,j,t}^{\textrm{con,tr}})$ is computed
by solving the following least-squares problem: \[(c_{i,j,t}^{\textrm{con,tr}},b_{i,j,t}^{\textrm{con,tr}})=\operatorname*{argmin}_{(\tilde{c}_{i,j,t},\tilde{b}_{i,j,t})}\sum_{k=1}^{p}\left(\tilde{c}_{i,j,t}(\bar{a}_{j}^{t}-\bar{d}_{i}^{t})^{(k)}+\tilde{b}_{i,j,t}-(\bar{E}_{i,j,t}^{\textrm{con,tr}})^{(k)}\right)^{2},\]
where $(i,j)\in\mathcal{A}^{t},\,t\in\mathcal{T}$; here and in what follows, bars distinguish observed data from decision variables, so $(\bar{E}_{i,j,t}^{\textrm{con,tr}})^{(k)}$ denotes the measured traction energy associated with the $k$-th observed track-based trip time $(\bar{a}_{j}^{t}-\bar{d}_{i}^{t})^{(k)}$, for $k=1,2,\ldots,p$. Likewise, $(c_{i,j,t,t^{\prime}}^{\textrm{con,cr}},b_{i,j,t,t^{\prime}}^{\textrm{con,cr}})$ is computed by solving the following least-squares problem:

\[
(c_{i,j,t,t^{\prime}}^{\textrm{con,cr}},b_{i,j,t,t^{\prime}}^{\textrm{con,cr}})=\operatorname*{argmin}_{(\tilde{c}_{i,j,t,t^{\prime}},\tilde{b}_{i,j,t,t^{\prime}})}\sum_{k=1}^{q}\left(\tilde{c}_{i,j,t,t^{\prime}}(\bar{a}_{j}^{t'}-\bar{d}_{i}^{t})^{(k)}+\tilde{b}_{i,j,t,t^{\prime}}-(\bar{E}_{i,j,t,t^{\prime}}^{\textrm{con,cr}})^{(k)}\right)^{2},\] where $(i,j)\in \varphi,(t,t^{\prime})\in\mathcal{B}_{ij}$ with $(\bar{E}_{i,j,t,t^{\prime}}^{\textrm{con,cr}})^{(k)}$ denoting the measured traction energy associated with the $k$-th observation of the crossing-over trip time $(\bar{a}_{j}^{t'}-\bar{d}_{i}^{t})^{(k)}$, for $k=1,2,\ldots,q$.

\subsection{Modeling transferred regenerative energy $E^{\textrm{reg}}$ \label{sec:regen_descrip}}

\begin{figure}
\begin{centering}
\includegraphics[scale=0.65]{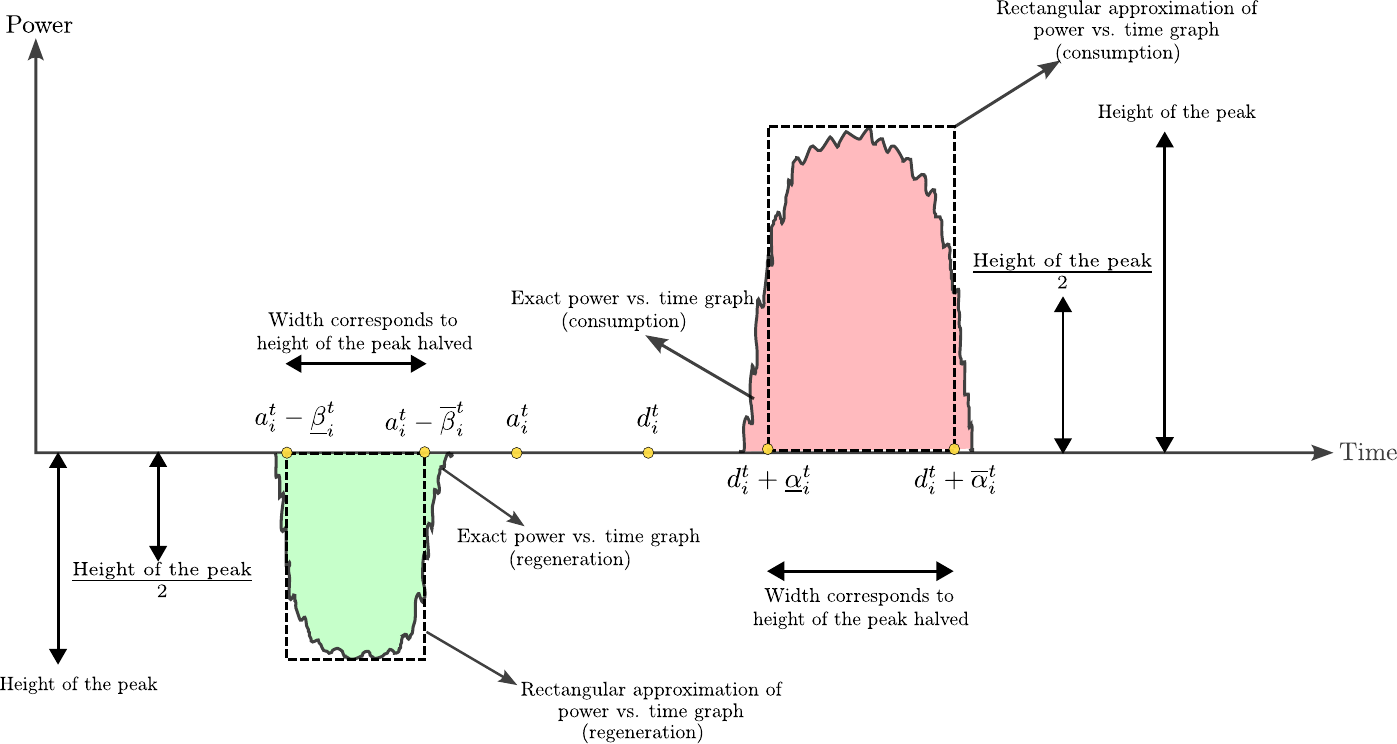}
\par\end{centering}
\caption{Effective braking and accelerating phases of a train $t$ at platform $i$, computed with the FWHM (Full Width at Half Maximum) rule. Each power profile is replaced by a rectangle whose height equals the peak power and whose width equals the full width at half maximum. The rectangles determine the endpoints $a_{i}^{t}-\underline{\beta}_{i}^{t}$ and $a_{i}^{t}-\overline{\beta}_{i}^{t}$ of the effective braking phase, and $d_{i}^{t}+\underline{\alpha}_{i}^{t}$ and $d_{i}^{t}+\overline{\alpha}_{i}^{t}$ of the effective accelerating phase. \label{Fig:FWHM}}
\end{figure}

\paragraph{Robust approximations of the power-versus-time graph.}
Maximizing regenerative-energy exchange over suitable train pairs is equivalent to maximizing the overlap between their consumption and regeneration power profiles. In CBTC systems, these power-versus-time curves depend on trip time and are highly nonlinear, with no convenient analytical form. Empirically, however, most of the power is concentrated near a peak, which motivates a robust lumped approximation via the FWHM (Full Width at Half Maximum) rule. This idea, popularized in early quantum mechanics \citep{gamow1985thirty} and now standard in areas such as numerical analysis, spectroscopy, signal processing, and statistics \citep[pp.~35--37]{Mahajan2010}, replaces each power profile by a rectangle whose height equals the peak value and whose width equals the full width at half maximum. In our setting, these rectangles are associated with the effective accelerating and braking intervals containing most of the energy as shown in Figure~\ref{Fig:FWHM}.

Consider a train $t$ arriving at (i.e., braking), then dwelling, and finally departing (i.e., accelerating) from platform $i$. As shown in Figure~\ref{Fig:FWHM}, the effective braking and accelerating phases of the train can be compactly contained in the FWHM rectangles. We denote the beginning and end of the effective braking phase of $t$ by $a_{i}^{t}-\underline{\beta}_{i}^{t}$ and $a_{i}^{t}-\overline{\beta}_{i}^{t}$, respectively, and the beginning and end of the effective accelerating phase of train $t$ by $d_{i}^{t}+\underline{\alpha}_{i}^{t}$ and $d_{i}^{t}+\overline{\alpha}_{i}^{t}$, respectively. Note the resulting ordering: $\underline{\beta}_{i}^{t}$ and $\overline{\beta}_{i}^{t}$ are both offsets measured backward from the arrival time, so the effective braking phase is the interval $[a_{i}^{t}-\underline{\beta}_{i}^{t},\,a_{i}^{t}-\overline{\beta}_{i}^{t}]$ with $\underline{\beta}_{i}^{t}\geq\overline{\beta}_{i}^{t}$, whereas the accelerating-phase offsets satisfy $\underline{\alpha}_{i}^{t}\leq\overline{\alpha}_{i}^{t}$.

Once these effective phase endpoints are known, the \textit{effective overlapping time} between two trains correspond to the duration for which the effective accelerating interval of one train overlaps the effective braking interval of the other, provided that the two platforms are supplied by the same electrical substation so that transmission losses can be neglected. Increasing this overlap increases the amount of regenerative energy that can be transferred between suitable train pairs, which we model next.

\paragraph{Constructing suitable train pairs for transfer of regenerative braking energy.}
We now identify train pairs that can realistically exchange regenerative energy. Let $\Omega$ collect the platform pairs between which such a transfer is physically possible: pairs supplied by the same electrical substation, so that transmission losses are negligible. For any $(i,j)\in\Omega$, let $\mathcal{T}_{i}\subseteq\mathcal{T}$ be the trains that arrive at, dwell at, and depart from platform $i$. To avoid duplicate platform pairs, we index the elements of $\Omega$ lexicographically with $i<j$. Fix $t\in\mathcal{T}_{i}$. We seek a train $\tilde{t}\in\mathcal{T}_{j}$ that, together with $t$, forms a suitable pair for regenerative-energy transfer. Our selection is based on an initial operating feasible timetable; in existing railway networks operated by CBTC systems such an operational timetable is always available. This initial feasible timetable provides feasible arrival and departure times $\bar{a}_{i}^{t}$ and $\bar{d}_{i}^{t}$ for every platform $i\in\mathcal{N}^{t}$; we write $\mu_{i}^{t}\triangleq(\bar{a}_{i}^{t}+\bar{d}_{i}^{t})/2$ for the midpoint of the corresponding dwell interval. Intuitively, among the trains that use platform $j$, the best candidate is the one whose dwell midpoint lies closest in time to that of $t$ in the reference timetable. This closeness can occur on either side of $t$, leading to the following definitions.

\begin{definition} \label{def:right-close} Fix $(i,j)\in\Omega$ and a train $t\in\mathcal{T}_{i}$. Consider the trains $s\in\mathcal{T}_{j}$ with $0\leq\mu_{j}^{s}-\mu_{i}^{t}\leq r$, where $r$ is an empirical parameter chosen by the timetable designer and much smaller than the full timetable horizon. If such trains exist, the one minimizing $\mu_{j}^{s}-\mu_{i}^{t}$, with ties broken by the trains' index order, is called \textit{temporally close to the right of $t$} and is denoted by $\overset{\rightharpoonup}{t}$; otherwise $t$ has no right event on $(i,j)$.
\end{definition}

\begin{definition} \label{def:left-close} Fix $(i,j)\in\Omega$ and a train $t\in\mathcal{T}_{i}$. Consider the trains $s\in\mathcal{T}_{j}$ with $0<\mu_{i}^{t}-\mu_{j}^{s}\leq r$. If such trains exist, the one minimizing $\mu_{i}^{t}-\mu_{j}^{s}$, with ties broken by the trains' index order, is called \textit{temporally close to the left of $t$} and is denoted by $\overset{\leftharpoonup}{t}$; otherwise $t$ has no left event on $(i,j)$. The left inequality is strict so that a pair of trains with coincident midpoints generates only a right event.
\end{definition}

\begin{definition} Fix $(i,j)\in\Omega$. For any train $t\in\mathcal{T}_{i}$, a train $\tilde{t}\in\mathcal{T}_{j}$ is \textit{temporally close to $t$} if it is temporally close to the right or to the left of $t$. \end{definition}

Thus, each potential synchronization event is identified by a tuple $(i,j,t,\tilde{t})$ satisfying the definitions above; because the platform pairs in $\Omega$ satisfy $i<j$, no event is enumerated twice. We collect these tuples in the set $\mathcal{E}$ and partition $\mathcal{E}$ into the subsets $\overset{\rightharpoonup}{\mathcal{E}}$ and $\overset{\leftharpoonup}{\mathcal{E}}$, which contain tuples of the forms $(i,j,t,\overset{\rightharpoonup}{t})$ and $(i,j,t,\overset{\leftharpoonup}{t})$, respectively. For a right event $(i,j,t,\overset{\rightharpoonup}{t})\in\overset{\rightharpoonup}{\mathcal{E}}$, we seek to align the effective accelerating phase of $t$ with the effective braking phase of $\overset{\rightharpoonup}{t}$. For a left event $(i,j,t,\overset{\leftharpoonup}{t})\in\overset{\leftharpoonup}{\mathcal{E}}$, the corresponding alignment is between the effective accelerating phase of $\overset{\leftharpoonup}{t}$ and the effective braking phase of $t$. The associated effective overlapping times are denoted by $\sigma_{ij}^{t\overset{\rightharpoonup}{t}}$ and $\sigma_{ij}^{t\overset{\leftharpoonup}{t}}$, respectively. An overlap may be positive when the two effective phases intersect, zero when they just meet, or negative when they are separated in time; these cases will be illustrated graphically later in this section.

The amount of regenerative energy that can be transferred depends on the total effective overlap, so our objective is to maximize that transfer. We write it as

\begin{align}
E^{\textrm{reg}}\triangleq\sum_{(i,j,t,\overset{\rightharpoonup}{t})\in\overset{\rightharpoonup}{\mathcal{E}}}\underbrace{E_{i,j,t,\overset{\rightharpoonup}{t}}^{\textrm{reg}}}_{\substack{\textrm{depends on }\\
\textrm{right event }\\
\textrm{overlapping time }\sigma_{ij}^{t\overset{\rightharpoonup}{t}}
}
}\quad+\sum_{(i,j,t,\overset{\leftharpoonup}{t})\in\overset{\leftharpoonup}{\mathcal{E}}}\underbrace{E_{i,j,t,\overset{\leftharpoonup}{t}}^{\textrm{reg}}}_{\substack{\textrm{depends on }\\
\textrm{left event }\\
\textrm{overlapping time }\sigma_{ij}^{t\overset{\leftharpoonup}{t}}
}
},
\label{eq:obj_funct-1}
\end{align}
where (i) $E_{i,j,t,\overset{\rightharpoonup}{t}}^{\textrm{reg}}:\mathbb{R}\to\mathbb{R}$ maps the overlap between the accelerating phase of $t$ and the braking phase of $\overset{\rightharpoonup}{t}$ to the corresponding transferred regenerative energy, with argument $\sigma_{ij}^{t\overset{\rightharpoonup}{t}}$, and (ii) $E_{i,j,t,\overset{\leftharpoonup}{t}}^{\textrm{reg}}:\mathbb{R}\to\mathbb{R}$ maps the overlap between the accelerating phase of $\overset{\leftharpoonup}{t}$ and the braking phase of $t$ to transferred regenerative energy, with argument $\sigma_{ij}^{t\overset{\leftharpoonup}{t}}$. To model \eqref{eq:obj_funct-1} tractably, we proceed in three steps: we fit the transferred-energy functions (Step 1), fit the effective accelerating and braking phases (Step 2), and express the overlaps in the decision variables (Step 3). We now describe each step in turn.

\paragraph{Step 1. Fitting the transferred-energy functions.}

Using a data-driven fit, we model $E_{i,j,t,\overset{\rightharpoonup}{t}}^{\textrm{reg}}$ and $E_{i,j,t,\overset{\leftharpoonup}{t}}^{\textrm{reg}}$ as affine functions of the effective overlapping times $\sigma_{ij}^{t\overset{\rightharpoonup}{t}}$ and $\sigma_{ij}^{t\overset{\leftharpoonup}{t}}$, respectively. This affine dependence is natural under the FWHM approximation: within the rectangles, power is nearly constant at its peak value, so the transferred energy scales with the overlap duration. Assumption~\ref{assum:triptime} keeps the overlaps within the narrow range where such fits are accurate. Thus, we write

\begin{equation}
\begin{alignedat}{1}E_{i,j,t,\overset{\rightharpoonup}{t}}^{\textrm{reg}} & \triangleq c_{i,j,t,\overset{\rightharpoonup}{t}}^{\textrm{reg}}\sigma_{ij}^{t\overset{\rightharpoonup}{t}}+b_{i,j,t,\overset{\rightharpoonup}{t}}^{\textrm{reg}},\\
E_{i,j,t,\overset{\leftharpoonup}{t}}^{\textrm{reg}} & \triangleq c_{i,j,t,\overset{\leftharpoonup}{t}}^{\textrm{reg}}\sigma_{ij}^{t\overset{\leftharpoonup}{t}}+b_{i,j,t,\overset{\leftharpoonup}{t}}^{\textrm{reg}},
\end{alignedat}
\label{eq:regen-energy-approx}
\end{equation}
where $(c_{i,j,t,\overset{\rightharpoonup}{t}}^{\textrm{reg}}, b_{i,j,t,\overset{\rightharpoonup}{t}}^{\textrm{reg}})$ and $(c_{i,j,t,\overset{\leftharpoonup}{t}}^{\textrm{reg}},b_{i,j,t,\overset{\leftharpoonup}{t}}^{\textrm{reg}})$ are nonnegative coefficients obtained from a least-squares fit constrained to a nonnegative slope and intercept, whose inputs are effective overlapping times and whose outputs are transferred regenerative-energy measurements. These input data may come either from historical operations or from physics-based simulation; in our numerical study they are generated by the \textsf{KDC} simulator. Positive intercepts $b_{i,j,t,\overset{\rightharpoonup}{t}}^{\textrm{reg}}$ and $b_{i,j,t,\overset{\leftharpoonup}{t}}^{\textrm{reg}}$ represent residual regenerative transfer, which can occur because some energy lies outside the FWHM rectangles even when the overlap is zero. Hence, \eqref{eq:regen-energy-approx} gives a practical interpretation of the objective in \eqref{eq:obj_funct-1}.

Since the fitted coefficients are nonnegative, positive overlap yields positive transferred regenerative energy. A nonpositive overlap, however, may still correspond to a small amount of transfer. When $\sigma_{ij}^{t\overset{\rightharpoonup}{t}}\leq-b_{i,j,t,\overset{\rightharpoonup}{t}}^{\textrm{reg}}/c_{i,j,t,\overset{\rightharpoonup}{t}}^{\textrm{reg}}$ (or $\sigma_{ij}^{t\overset{\leftharpoonup}{t}}\leq-b_{i,j,t,\overset{\leftharpoonup}{t}}^{\textrm{reg}}/c_{i,j,t,\overset{\leftharpoonup}{t}}^{\textrm{reg}}$), the transfer is negligible or zero. In those cases, the model allows negative values of $E_{i,j,t,\overset{\rightharpoonup}{t}}^{\textrm{reg}}$ (or $E_{i,j,t,\overset{\leftharpoonup}{t}}^{\textrm{reg}}$) to act as a penalty in the maximization objective. In our model, this penalty is an intentional modeling choice in contrast with truncating the transfer at zero inside the model because, even when the overlap is negative, we still want the separation as small as possible, so that some residual synchronization can be obtained through the optimization process.

\paragraph{Step 2. Fitting the effective accelerating and braking phases.}

To evaluate effective overlap, we need the start and end times of the effective accelerating and braking intervals. Consider train $t\in\mathcal{T}$ traveling from platform $i$ to platform $j$ on track $(i,j)\in\mathcal{A}^{t}$. The effective accelerating phase after departure from platform $i$ begins at $d_{i}^{t}+\underline{\alpha}_{i}^{t}$ and ends at $d_{i}^{t}+\overline{\alpha}_{i}^{t}$. The effective braking phase before arrival at platform $j$ begins at $a_{j}^{t}-\underline{\beta}_{j}^{t}$ and ends at $a_{j}^{t}-\overline{\beta}_{j}^{t}$. The quantities $\underline{\alpha}_{i}^{t},\overline{\alpha}_{i}^{t},\underline{\beta}_{j}^{t},\overline{\beta}_{j}^{t}$ depend on the train's power profile, and hence on the trip time $a_{j}^{t}-d_{i}^{t}$. Exact analytical expressions for these quantities as functions of trip time are not available, but Assumption~\ref{assum:triptime} again motivates affine approximations. Therefore, for $t\in\mathcal{T}$ and $(i,j)\in\mathcal{A}^{t}$, we define $\underline{\alpha}_{i}^{t},\overline{\alpha}_{i}^{t},\underline{\beta}_{j}^{t},\overline{\beta}_{j}^{t}$ as follows:
\begin{equation}
\begin{aligned}\underline{\alpha}_{i}^{t} & \triangleq c_{i,t}^{\underline{\alpha}}(a_{j}^{t}-d_{i}^{t})+b_{i,t}^{\underline{\alpha}},\quad\overline{\alpha}_{i}^{t} \triangleq c_{i,t}^{\overline{\alpha}}(a_{j}^{t}-d_{i}^{t})+b_{i,t}^{\overline{\alpha}},\\
\underline{\beta}_{j}^{t} & \triangleq c_{i,t}^{\underline{\beta}}(a_{j}^{t}-d_{i}^{t})+b_{i,t}^{\underline{\beta}},\quad\overline{\beta}_{j}^{t} \triangleq c_{i,t}^{\overline{\beta}}(a_{j}^{t}-d_{i}^{t})+b_{i,t}^{\overline{\beta}},
\end{aligned} \tag{\texttt{EFF\_ACCEL\_\&\_BRK\_PHASES}}
\label{eq:effective_braking_accelearation_compute}
\end{equation}
where the coefficient pairs $(c_{i,t}^{\underline{\alpha}},b_{i,t}^{\underline{\alpha}})$, $(c_{i,t}^{\overline{\alpha}},b_{i,t}^{\overline{\alpha}})$, $(c_{i,t}^{\underline{\beta}},b_{i,t}^{\underline{\beta}})$, and $(c_{i,t}^{\overline{\beta}},b_{i,t}^{\overline{\beta}})$ are obtained by least squares, using trip time as the input measurement and the corresponding phase start/end times as the output measurements. For each train $t$ and platform $j$, the incoming track $(i,j)\in\mathcal{A}^{t}$ is unique, so indexing these coefficients by $(i,t)$ is unambiguous.

\begin{figure}
\begin{centering}
\includegraphics[scale=0.58]{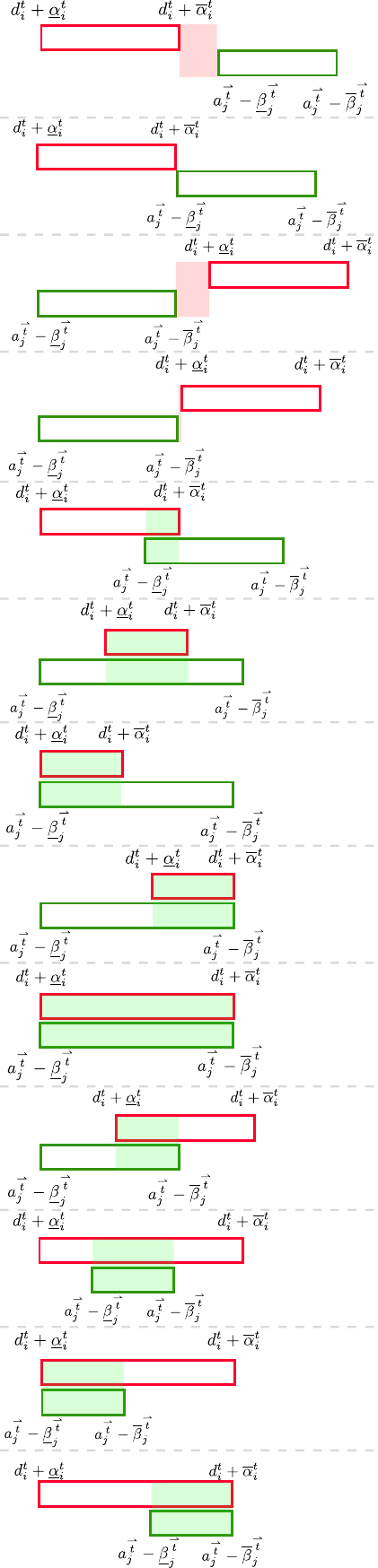}
\par\end{centering}
\caption{The thirteen possible overlapping times $\sigma_{ij}^{t\overset{\rightharpoonup}{t}}$ between the accelerating phase of train $t$ and the braking phase of train $\overset{\rightharpoonup}{t}$. In the first four cases, $\sigma_{ij}^{t\overset{\rightharpoonup}{t}}$ is nonpositive (red-shaded region); in the remaining nine cases, it is positive (green-shaded region). All thirteen cases obey the closed form $\sigma_{ij}^{t\overset{\rightharpoonup}{t}}=\min\left\{ a_{j}^{\overset{\rightharpoonup}{t}}-\overline{\beta}_{j}^{\overset{\rightharpoonup}{t}},d_{i}^{t}+\overline{\alpha}_{i}^{t}\right\} +\min\left\{ -d_{i}^{t}-\underline{\alpha}_{i}^{t},-a_{j}^{\overset{\rightharpoonup}{t}}+\underline{\beta}_{j}^{\overset{\rightharpoonup}{t}}\right\}$. \label{fig:overlap}}
\end{figure}

\paragraph{Step 3. Expressing the overlaps in the decision variables.}

Now, we describe how to model the effective overlapping times $\sigma_{ij}^{t\overset{\rightharpoonup}{t}}$ and $\sigma_{ij}^{t\overset{\leftharpoonup}{t}}$ in terms of the decision variables.

First, we discuss the modeling of the right-event overlapping time $\sigma_{ij}^{t\overset{\rightharpoonup}{t}}$.
By Allen's interval algebra \citep{Allen1983}, exactly thirteen temporal relations
are possible between the accelerating phase
of train $t$ and the braking phase of train $\overset{\rightharpoonup}{t}$;
they are shown in Figure \ref{fig:overlap}. In the first four cases, the overlap is nonpositive: the two intervals are disjoint (red-shaded); in the remaining nine, it is positive (green-shaded). Fortunately, we can write
all thirteen possible overlapping times in closed form:
\begin{align}
\sigma_{ij}^{t\overset{\rightharpoonup}{t}} & =\min\left\{ a_{j}^{\overset{\rightharpoonup}{t}}-\overline{\beta}_{j}^{\overset{\rightharpoonup}{t}},d_{i}^{t}+\overline{\alpha}_{i}^{t}\right\} -\max\left\{ d_{i}^{t}+\underline{\alpha}_{i}^{t},a_{j}^{\overset{\rightharpoonup}{t}}-\underline{\beta}_{j}^{\overset{\rightharpoonup}{t}}\right\} \nonumber \\
 & =\min\left\{ a_{j}^{\overset{\rightharpoonup}{t}}-\overline{\beta}_{j}^{\overset{\rightharpoonup}{t}},d_{i}^{t}+\overline{\alpha}_{i}^{t}\right\} +\min\left\{ -d_{i}^{t}-\underline{\alpha}_{i}^{t},-a_{j}^{\overset{\rightharpoonup}{t}}+\underline{\beta}_{j}^{\overset{\rightharpoonup}{t}}\right\} ,\label{eq:sigma_right} \tag{\texttt{RGHT\_EVNT\_OVLP\_TIME}}
\end{align}
which is concave in the decision variables, being a pointwise minimum of affine functions \citep[Section~3.2]{Boyd2004}. For example, suppose train $t$ departs platform $i$ at $d_{i}^{t}=100$ s with $\underline{\alpha}_{i}^{t}=4$ s and $\overline{\alpha}_{i}^{t}=18$ s, while train $\overset{\rightharpoonup}{t}$ arrives at platform $j$ at $a_{j}^{\overset{\rightharpoonup}{t}}=120$ s with $\underline{\beta}_{j}^{\overset{\rightharpoonup}{t}}=15$ s and $\overline{\beta}_{j}^{\overset{\rightharpoonup}{t}}=3$ s. Then the effective accelerating phase is $[104,118]$, the effective braking phase is $[105,117]$, and \eqref{eq:sigma_right} gives $\sigma_{ij}^{t\overset{\rightharpoonup}{t}}=\min\{117,118\}-\max\{104,105\}=12$ s, exactly the length of time the two phases share. Next, we discuss the modeling of the left-event overlapping time $\sigma_{ij}^{t\overset{\leftharpoonup}{t}}$. Replacing $i$, $j$, $t$, and $\overset{\rightharpoonup}{t}$
in \eqref{eq:sigma_right} with $j$, $i$,
$\overset{\leftharpoonup}{t}$, and $t$, respectively, gives
\begin{align}
\sigma_{ij}^{t\overset{\leftharpoonup}{t}} & =\min\left\{ a_{i}^{t}-\overline{\beta}_{i}^{t},d_{j}^{\overset{\leftharpoonup}{t}}+\overline{\alpha}_{j}^{\overset{\leftharpoonup}{t}}\right\} -\max\left\{ d_{j}^{\overset{\leftharpoonup}{t}}+\underline{\alpha}_{j}^{\overset{\leftharpoonup}{t}},a_{i}^{t}-\underline{\beta}_{i}^{t}\right\} \nonumber \\
 & =\min\left\{ a_{i}^{t}-\overline{\beta}_{i}^{t},d_{j}^{\overset{\leftharpoonup}{t}}+\overline{\alpha}_{j}^{\overset{\leftharpoonup}{t}}\right\} +\min\left\{ -d_{j}^{\overset{\leftharpoonup}{t}}-\underline{\alpha}_{j}^{\overset{\leftharpoonup}{t}},-a_{i}^{t}+\underline{\beta}_{i}^{t}\right\} ,\label{eq:sigma_left}  \tag{\texttt{LEFT\_EVNT\_OVLP\_TIME}}
\end{align}
which is again a concave function of our decision variables. Because each overlap is concave in the event times and each Step~1 fit is affine with
a nonnegative slope, the total transferred regenerative energy is a concave function
of the timetable; Section~\ref{sec:Final-optimization-model} maximizes it through
linear hypograph constraints.

\section{Final optimization model \label{sec:Final-optimization-model}}

Combining the ingredients developed above yields the following linear optimization model for minimizing the railway network's effective energy consumption:
\begin{equation}
\begin{aligned}& \underset{}{\mbox{minimize}}\\
& \sum_{t\in\mathcal{T},(i,j)\in\mathcal{A}^{t}}\left(c_{i,j,t}^{\textrm{con,tr}}(a_{j}^{t}-d_{i}^{t})+b_{i,j,t}^{\textrm{con,tr}}\right)+\sum_{(i,j)\in\varphi,(t,t^{\prime})\in\mathcal{B}_{ij}}\left(c_{i,j,t,t^{\prime}}^{\textrm{con,cr}}(a_{j}^{t^{\prime}}-d_{i}^{t})+b_{i,j,t,t^{\prime}}^{\textrm{con,cr}}\right) \;\;\rhd\;\textrm{consumption}\\
& -\sum_{(i,j,t,\overset{\rightharpoonup}{t})\in\overset{\rightharpoonup}{\mathcal{E}}}\left(c_{i,j,t,\overset{\rightharpoonup}{t}}^{\textrm{reg}}\sigma_{ij}^{t\overset{\rightharpoonup}{t}}+b_{i,j,t,\overset{\rightharpoonup}{t}}^{\textrm{reg}}\right)-\sum_{(i,j,t,\overset{\leftharpoonup}{t})\in\overset{\leftharpoonup}{\mathcal{E}}}\left(c_{i,j,t,\overset{\leftharpoonup}{t}}^{\textrm{reg}}\sigma_{ij}^{t\overset{\leftharpoonup}{t}}+b_{i,j,t,\overset{\leftharpoonup}{t}}^{\textrm{reg}}\right) \;\rhd\;\textrm{transferred regen. energy}\\
& \textup{subject to} \\
& \begin{rcases}
\sigma_{ij}^{t\overset{\rightharpoonup}{t}}\leq\varpi_{ij}^{t\overset{\rightharpoonup}{t}}+\psi_{ij}^{t\overset{\rightharpoonup}{t}}\\
\varpi_{ij}^{t\overset{\rightharpoonup}{t}}\leq a_{j}^{\overset{\rightharpoonup}{t}}-\overline{\beta}_{j}^{\overset{\rightharpoonup}{t}}\\
\varpi_{ij}^{t\overset{\rightharpoonup}{t}}\leq d_{i}^{t}+\overline{\alpha}_{i}^{t}\\
\psi_{ij}^{t\overset{\rightharpoonup}{t}}\leq-d_{i}^{t}-\underline{\alpha}_{i}^{t}\\
\psi_{ij}^{t\overset{\rightharpoonup}{t}}\leq-a_{j}^{\overset{\rightharpoonup}{t}}+\underline{\beta}_{j}^{\overset{\rightharpoonup}{t}}
\end{rcases} \textrm{ for }(i,j,t,\overset{\rightharpoonup}{t})\in\overset{\rightharpoonup}{\mathcal{E}} \quad \substack{\rhd \; \textrm{hypograph form}\\
\textrm{of }\eqref{eq:sigma_right}
} \\
& \begin{rcases}
\sigma_{ij}^{t\overset{\leftharpoonup}{t}}\leq\varpi_{ij}^{t\overset{\leftharpoonup}{t}}+\psi_{ij}^{t\overset{\leftharpoonup}{t}}\\
\varpi_{ij}^{t\overset{\leftharpoonup}{t}}\leq a_{i}^{t}-\overline{\beta}_{i}^{t}\\
\varpi_{ij}^{t\overset{\leftharpoonup}{t}}\leq d_{j}^{\overset{\leftharpoonup}{t}}+\overline{\alpha}_{j}^{\overset{\leftharpoonup}{t}}\\
\psi_{ij}^{t\overset{\leftharpoonup}{t}}\leq-d_{j}^{\overset{\leftharpoonup}{t}}-\underline{\alpha}_{j}^{\overset{\leftharpoonup}{t}}\\
\psi_{ij}^{t\overset{\leftharpoonup}{t}}\leq-a_{i}^{t}+\underline{\beta}_{i}^{t}\\
\end{rcases} \textrm{ for }(i,j,t,\overset{\leftharpoonup}{t})\in\overset{\leftharpoonup}{\mathcal{E}} \quad \substack{\rhd \; \textrm{hypograph form}\\
\textrm{of }\eqref{eq:sigma_left}
} \\
& \textrm{Eq.} \; \eqref{eq:effective_braking_accelearation_compute}, \; \rhd \; \textrm{effective-phase endpoints} \\
& \textrm{Eq.}\; \eqref{eq:TripTimeConstraint}, \eqref{eq:turnAroundConstraint}, \eqref{eq:DwellTimeConstraint}, \eqref{eq:ConnectionConstraint}, \; \rhd \; \textrm{operational constraints} \\
&  \textrm{Eq.}\;  \eqref{eq:SafetyConstraint1}, \eqref{eq:TotalTravelTimeConstraints}, \eqref{eq:domain}, \; \rhd \; \textrm{headway, total travel time, and domain}
\end{aligned} \tag{\texttt{OPT\_MODEL}} \label{eq:final_lp_model}
\end{equation}
where the decision variables are $a_{i}^{t}$, $d_{i}^{t}$, $\sigma_{ij}^{t\overset{\rightharpoonup}{t}},$$\varpi_{ij}^{t\overset{\rightharpoonup}{t}},$
$\psi_{ij}^{t\overset{\rightharpoonup}{t}}$, $\sigma_{ij}^{t\overset{\leftharpoonup}{t}},$
$\varpi_{ij}^{t\overset{\leftharpoonup}{t}},$ and $\psi_{ij}^{t\overset{\leftharpoonup}{t}}$. Note that the first ten constraints recast \eqref{eq:sigma_right} and \eqref{eq:sigma_left} as linear constraints using the hypograph approach \citep[Section~4.1.3]{Boyd2004}. Because the objective coefficients $c^{\textrm{reg}}$ are nonnegative, the optimization drives each auxiliary variable up to its bound. Hence, at an optimal solution, every overlap variable with a positive coefficient equals its overlap value in \eqref{eq:sigma_right} or \eqref{eq:sigma_left}. In words, \eqref{eq:final_lp_model} chooses the arrival and departure times, subject to all the operational constraints, to make the predicted traction energy minus the transferred regenerative energy as small as possible.
 
 \paragraph{Predicting effective energy consumption from the solution.} After solving \eqref{eq:final_lp_model}, we obtain an energy-optimal timetable. Its predicted effective energy consumption is

\begin{align}
 & \sum_{t\in\mathcal{T},(i,j)\in\mathcal{A}^{t}}\left(c_{i,j,t}^{\textrm{con,tr}}(a_{j}^{t}-d_{i}^{t})+b_{i,j,t}^{\textrm{con,tr}}\right)\;\quad+\sum_{(i,j)\in\varphi,(t,t^{\prime})\in\mathcal{B}_{ij}}\left(c_{i,j,t,t^{\prime}}^{\textrm{con,cr}}(a_{j}^{t^{\prime}}-d_{i}^{t})+b_{i,j,t,t^{\prime}}^{\textrm{con,cr}}\right)\nonumber \\
 & -\sum_{(i,j,t,\overset{\rightharpoonup}{t})\in\overset{\rightharpoonup}{\mathcal{E}}}\max\left\{ c_{i,j,t,\overset{\rightharpoonup}{t}}^{\textrm{reg}}\sigma_{ij}^{t\overset{\rightharpoonup}{t}}+b_{i,j,t,\overset{\rightharpoonup}{t}}^{\textrm{reg}},0\right\} \;\quad-\sum_{(i,j,t,\overset{\leftharpoonup}{t})\in\overset{\leftharpoonup}{\mathcal{E}}}\max\left\{ c_{i,j,t,\overset{\leftharpoonup}{t}}^{\textrm{reg}}\sigma_{ij}^{t\overset{\leftharpoonup}{t}}+b_{i,j,t,\overset{\leftharpoonup}{t}}^{\textrm{reg}},0\right\} .\tag{\texttt{PRED\_EFF\_ENG\_CNSM}}\label{eq:predicted_eng}
\end{align}
In \eqref{eq:predicted_eng}, the first line represents traction energy consumed during acceleration, while the second line subtracts the regenerative energy that is actually transferred. The full expression therefore gives the predicted effective energy consumption. When $\sigma_{ij}^{t\overset{\rightharpoonup}{t}}\leq0$ and $\sigma_{ij}^{t\overset{\leftharpoonup}{t}}\leq0$, the affine terms $c_{i,j,t,\overset{\rightharpoonup}{t}}^{\textrm{reg}}\sigma_{ij}^{t\overset{\rightharpoonup}{t}}+b_{i,j,t,\overset{\rightharpoonup}{t}}^{\textrm{reg}}$ and $c_{i,j,t,\overset{\leftharpoonup}{t}}^{\textrm{reg}}\sigma_{ij}^{t\overset{\leftharpoonup}{t}}+b_{i,j,t,\overset{\leftharpoonup}{t}}^{\textrm{reg}}$ act as penalties in the objective of \eqref{eq:final_lp_model}; accordingly, the clipped quantities $\max\{c_{i,j,t,\overset{\rightharpoonup}{t}}^{\textrm{reg}}\sigma_{ij}^{t\overset{\rightharpoonup}{t}}+b_{i,j,t,\overset{\rightharpoonup}{t}}^{\textrm{reg}},0\}$ and $\max\{c_{i,j,t,\overset{\leftharpoonup}{t}}^{\textrm{reg}}\sigma_{ij}^{t\overset{\leftharpoonup}{t}}+b_{i,j,t,\overset{\leftharpoonup}{t}}^{\textrm{reg}},0\}$ recover the actual regenerative energy transferred in the individual right and left events, respectively. Moreover, each affine transfer term is bounded above by its clipped counterpart, so the objective value of \eqref{eq:final_lp_model} is an upper bound on \eqref{eq:predicted_eng}: minimizing the former minimizes an upper bound on the predicted effective energy consumption.

\section{Numerical experiments \label{Numerical_Study} }

We evaluate the model on two real-world railway networks. Section~\ref{sec:sl8} describes Metro Line 8 of the Shanghai Railway Network and presents the corresponding numerical results; Section~\ref{sec:dlr} applies the same analysis to the Docklands Light Railway.

We ran all experiments in this section on a 16-inch MacBook Pro with an 18-core Apple M5 Max processor, 128~GB of unified memory, and macOS 26.5.2. We used Julia 1.12.6, \textsf{JuMP} 1.30.1 \citep{Lubin2023}, and \textsf{Gurobi Optimizer} 13.0.2 \citep{gurobi2026} for solving each linear program. The solution times in Tables~\ref{tab:Results-of-the-numerical-study-sl8} and~\ref{tab:Results-of-the-numerical-study-dlr} are the medians of three Gurobi-reported solver times obtained from three independent runs; each minimization solve was warm-started at the corresponding baseline timetable. The synchronization events, the platform-pair set $\Omega$, and the fitted energy coefficients are supplied with the instance data. In the Shanghai instances, the synchronization events were constructed with the closeness parameter $r=120$~seconds; in the Docklands Light Railway instances, the paired dwell midpoints lie within 10~minutes of each other in the reference timetables. In every instance of both networks, no braking train's effective phase enters more than one synchronization event, so each transfer draws on a distinct braking phase.

\begin{figure}
\begin{centering}
\includegraphics[scale=0.45]{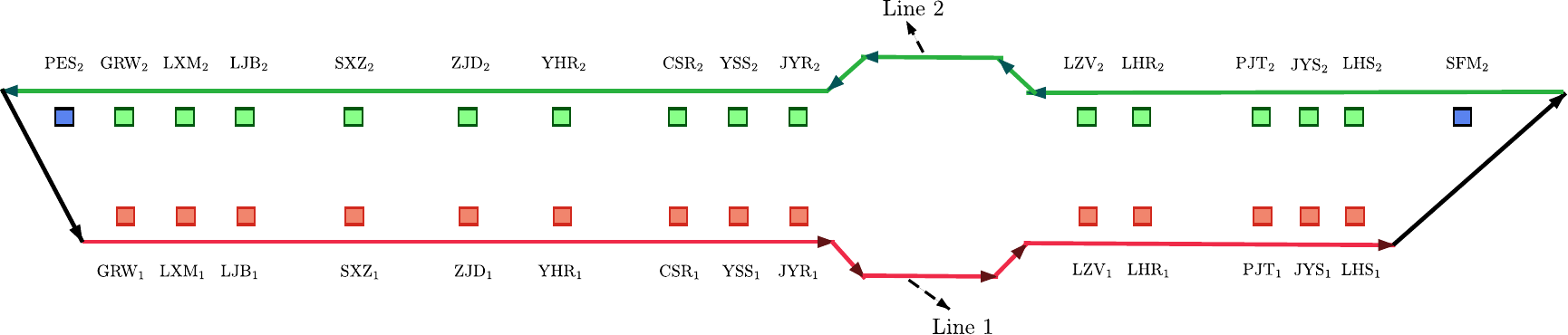}
\par\end{centering}
\caption{Shanghai Metro Line 8 network used in the numerical study. Squares denote platforms: red on Line~1, green on Line~2; the blue squares are the terminal platforms joined to the opposite line by crossing-overs.\label{Fig:sl8}}
\end{figure}

\subsection{Metro Line 8 in Shanghai Railway Network}\label{sec:sl8}

\subsubsection{Railway network in consideration}\label{sec:NetworkSl8}

We run our experiments on the $\textrm{PES}_2$-$\textrm{SFM}_2$ service of Metro Line 8 in the Shanghai Railway Network, which is the world's longest metro system by route length and carried approximately 3.88 billion passenger trips in 2019 \citep{UITP2022}. Metro Line 8 serves some of Shanghai's densest residential districts and carried an average daily ridership of about 1.08 million at the time of the study of \citet{gupta2016two}. The line extends 37.4~km and has 28 operating stations \citep{ShmetroLine8}.

Three services run on Metro Line 8, each with its own CBTC-generated timetable; here we optimize the $\textrm{PES}_2$-$\textrm{SFM}_2$ service shown in Figure~\ref{Fig:sl8}. The corresponding network consists of two lines, denoted Line 1 and Line 2, and 14 stations, labeled in uppercase in the figure. Each station has one platform on each line; for example, LXM includes platforms $\textrm{LXM}_1$ and $\textrm{LXM}_2$ serving Line 1 and Line 2, respectively. The portion of the line considered here is 37 km long. The mean inter-station spacing is 1.4 km, with a minimum of 738 m between YHR and ZJD and a maximum of 2.6 km between PJT and LHR. Track grades lie in $[-2.00453^{\circ},2.00453^{\circ}]$, the maximum acceleration is 1.04~$\textrm{m/s}^{2}$, and the maximum braking deceleration is $0.8~\textrm{m/s}^{2}$. The electrical-to-kinetic conversion efficiency is 0.9, the kinetic-to-regenerative-braking conversion efficiency is 0.76, and the regenerative-electricity transmission-loss factor is 0.1. 

The affine energy models underlying these instances were fitted offline from \textsf{KDC} output; generating the training data took roughly one hour of simulation per network, and all least-squares fits complete in under a minute. 
The 29 track-specific traction fits are shared by all trains traversing the corresponding track and by all 11 Shanghai instances. Their coefficients remain fixed; the traction energy predicted for a traversal depends on its timetable-dependent trip time. Five simulated trip times per track, chosen across the admissible window,
give coefficients of determination with median $0.98$ (fifth percentile $0.83$, minimum $0.76$) and relative root-mean-square errors of at most $4.2\%$ over the 29 track fits (the crossing-over fits behave similarly, with coefficients of determination between $0.82$ and $0.94$). The robust trip-time windows are themselves narrow: $0.4$--$21.2$~seconds wide for \eqref{eq:TripTimeConstraint} (median $6.8$) and $8.4$--$31.4$~seconds for \eqref{eq:turnAroundConstraint}. This makes Assumption~\ref{assum:triptime} quantitative for these instances. The 25 transfer fits use 42--108 measurements each, covering effective overlaps from $-32.9$ to $23.0$~seconds; their coefficients of determination have median $0.90$ (fifth percentile $0.70$, minimum $0.32$), all fitted slopes and intercepts are nonnegative, and the aggregate effect of the residual fitting error on the reported reductions is bounded by the independent \textsf{SPSIM} evaluation of Section~\ref{sec:resultsSl8}. The fitted slopes of the effective-phase offsets in \eqref{eq:effective_braking_accelearation_compute} are negligible over windows this narrow, so the instance data supply the offsets as constants, and the experiments in effect use constant effective phases.

\subsubsection{Results of the numerical study}\label{sec:resultsSl8}

We evaluate the model on 11 operational instances of the $\textrm{PES}_2$-$\textrm{SFM}_2$ service on Metro Line 8. The instances correspond to the eleven service configurations studied in \citet{gupta2016two}, with baseline timetables and operational data from a more recent generation of the CBTC system. The number of trains and headway durations vary across instances. Hitachi Rail Canada supplied the corresponding baseline timetables, and we measure the benefit of the optimized schedules by comparing their energy usage against these existing timetables.

\begin{table*}[htbp]
  \centering
  \scriptsize 
  \caption{Summary of numerical results for Metro Line 8 of the Shanghai Railway Network. }
  \label{tab:Results-of-the-numerical-study-sl8}
  \setlength{\tabcolsep}{2.5pt}
  \begin{tabular}{C{1.7cm}C{1.7cm}C{2.05cm}C{1.9cm}C{1.9cm}C{1.9cm}C{1.9cm}C{1.9cm}}
    \toprule
    \textbf{\#Trains} & \textbf{\#Variables} & \textbf{\#Constraints} & \textbf{Solution time (s)} & \textbf{Initial effective energy consumption (kWh)} & \textbf{ Predicted final effective energy consumption (kWh)} & \textbf{Reduction predicted by  our model (\%)}  & \textbf{Reduction predicted by  \textsf{SPSIM} (\%)} \\
    \midrule
    1000 & 47581 & 151259 & 0.448 & 353721.62 & 254062.35 & 28.17 & 30.27 \\
    1032 & 49104 & 156102 & 0.432 & 364832.93 & 266932.98 & 26.83 & 31.34 \\
    1066 & 50726 & 161254 & 0.481 & 376017.98 & 274746.55 & 26.93 & 33.04\\
    1100 & 52339 & 166391 & 0.476 & 386087.42 & 288773.33 & 25.21 & 29.10\\
    1132 & 53865 & 171239 & 0.526 & 389844.65 & 294141.71 & 24.55 & 26.56\\
    1166 & 55487 & 176391 & 0.530 & 407717.20 & 307967.10 & 24.47 & 26.55\\
    1198 & 57010 & 181234 & 0.579 & 417767.32 & 336232.08 & 19.52 & 19.27\\
    1232 & 58626 & 186376 & 0.585 & 447025.23 & 340027.49 & 23.94 & 25.86\\
    1266 & 60242 & 191518 & 0.596 & 438657.08 & 337269.06 & 23.11 & 26.94 \\
    1298 & 61765 & 196361 & 0.651 & 449379.24 & 336794.26 & 25.05 & 24.83\\
    1332 & 63387 & 201513 & 0.665 & 475354.59 & 354114.70 & 25.51 & 25.97\\
    \bottomrule
  \end{tabular}
\end{table*}

\begin{figure}[ht]
    \centering
    \begin{subfigure}[b]{0.4\textwidth}
        \centering
        \includegraphics[width=\textwidth]{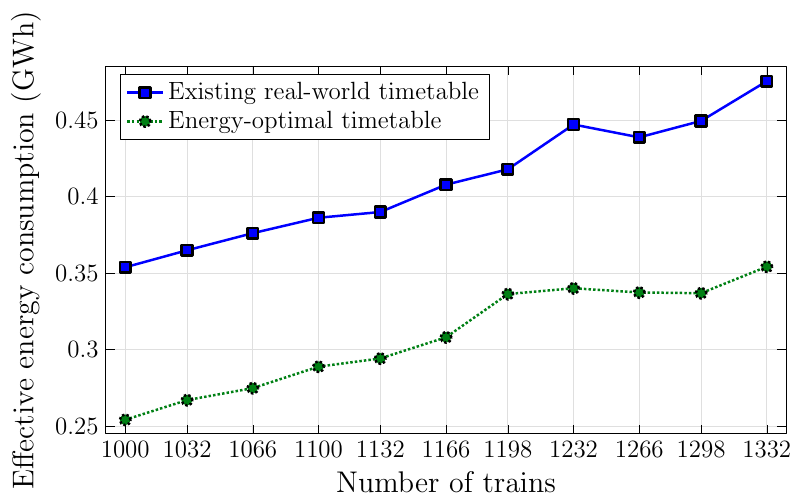}
        \caption{}
        \label{fig:a}
    \end{subfigure}
    \hfill
    \begin{subfigure}[b]{0.4\textwidth}
        \centering
        \includegraphics[width=\textwidth]{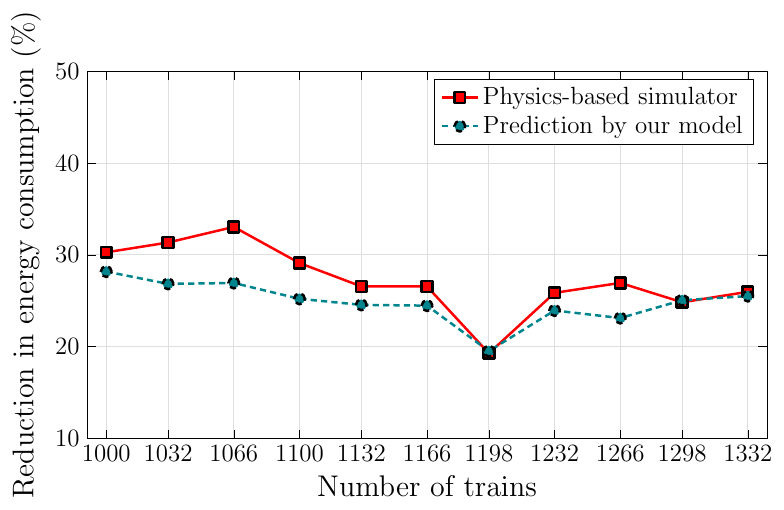}
        \caption{}
        \label{fig:b}
    \end{subfigure}
        \caption{Shanghai Metro Line 8 results. (a)~Model-predicted effective energy consumption of the existing and optimized timetables. (b)~Percentage reductions in effective energy consumption predicted by our model and computed with \textsf{SPSIM}. Each reduction compares the optimized timetable with the existing timetable using the same evaluator.}
    \label{fig:ab}
\end{figure}

The findings of this computational study are summarized in Table \ref{tab:Results-of-the-numerical-study-sl8}. Across the 11 instances, all minimization solves terminated optimally, and the median solution times range from 0.432 to 0.665 seconds.
Both the initial and optimized energy values in Table~\ref{tab:Results-of-the-numerical-study-sl8} and Figure~\ref{fig:a} are computed with \eqref{eq:predicted_eng}, evaluated at the corresponding timetable's event times. The model-predicted percentage reduction compares these two values using the model-evaluated baseline energy as its denominator. On the other hand, the \textsf{SPSIM} percentage reduction is computed independently by evaluating both timetables with the \textsf{SPSIM} simulator and using its baseline energy as the denominator.
Here we note that \textsf{SPSIM} and \textsf{KDC} simulators operate independently. \textsf{SPSIM} computes the network-level energy consumption of an entire timetable, whereas \textsf{KDC} (that we used in Section~\ref{sec:Modeling-the-objective}) generates the single-train power-profile data underlying our affine fits. After solving \eqref{eq:final_lp_model}, we evaluate the effective energy consumption of each optimized timetable with \eqref{eq:predicted_eng}. The optimized timetables improve substantially on their real-world counterparts, with predicted reductions in effective energy consumption ranging from 19.52\% in the least favorable instance to 28.17\% in the most favorable one. Figure~\ref{fig:b} plots these reductions.

To check how accurate the model's predictions are, we input each optimized timetable into \textsf{SPSIM} and let the simulator compute the effective energy consumption and the resulting reduction. The last two columns of Table~\ref{tab:Results-of-the-numerical-study-sl8} and Figure~\ref{fig:b} show that the two predictions agree closely, with our model being slightly more conservative on average. The gap between our model and \textsf{SPSIM} has a simple explanation: \textsf{SPSIM} computes transferred regenerative energy from the exact area of overlap between the power-versus-time curves of an accelerating--braking pair, whereas our model uses the rectangular approximations of Section~\ref{sec:regen_descrip} that sometimes capture less energy than the exact overlap does.

Table~\ref{tab:decomposition-sl8} decomposes the predicted reduction into its two mechanisms: (i) lower traction energy from trip times lengthened within their admissible windows, and (ii) additional regenerative transfer from synchronized braking and acceleration. Both mechanisms contribute separately in every instance: the traction component supplies 51.0\%--73.4\% of the reduction and the transfer component the remaining 26.6\%--49.0\%, with transferred regenerative energy growing by a factor of 1.6--3.4 over the baseline timetables.

\begin{table*}[htbp]
  \centering
  \scriptsize
  \caption{Decomposition of the predicted effective-energy reduction for Metro Line 8: traction energy $E^{\textrm{con}}$ and transferred regenerative energy $E^{\textrm{reg}}$ at the baseline and optimized timetables (kWh, evaluated with the model), and each mechanism's share of the reduction.}
  \label{tab:decomposition-sl8}
  \begin{tabular}{r rr rr rr}
    \toprule
    & \multicolumn{2}{c}{Baseline} & \multicolumn{2}{c}{Optimized} & \multicolumn{2}{c}{Share of reduction (\%)} \\
    \cmidrule(lr){2-3}\cmidrule(lr){4-5}\cmidrule(lr){6-7}
    \#Trains & $E^{\textrm{con}}$ & $E^{\textrm{reg}}$ & $E^{\textrm{con}}$ & $E^{\textrm{reg}}$ & traction & transfer \\
    \midrule
    1000 & 378545.86 & 24824.24 & 327207.60 & 73152.61 & 51.5 & 48.5 \\
    1032 & 390658.93 & 25826.00 & 337978.83 & 71045.85 & 53.8 & 46.2 \\
    1066 & 403529.69 & 27511.71 & 350675.36 & 75941.72 & 52.2 & 47.8 \\
    1100 & 416400.45 & 30313.03 & 361444.27 & 72669.77 & 56.5 & 43.5 \\
    1132 & 428513.51 & 38668.86 & 373276.48 & 79135.21 & 57.7 & 42.3 \\
    1166 & 441384.27 & 33667.08 & 383841.65 & 75874.71 & 57.7 & 42.3 \\
    1198 & 453497.34 & 35730.02 & 393653.84 & 57403.06 & 73.4 & 26.6 \\
    1232 & 466368.10 & 19342.87 & 405924.57 & 65833.13 & 56.5 & 43.5 \\
    1266 & 479238.86 & 40581.78 & 418401.17 & 81128.18 & 60.0 & 40.0 \\
    1298 & 491351.92 & 41972.68 & 427633.94 & 90795.08 & 56.6 & 43.4 \\
    1332 & 504222.69 & 28868.10 & 442391.52 & 88264.42 & 51.0 & 49.0 \\
    \bottomrule
  \end{tabular}
\end{table*}

\subsection{Docklands Light Railway}\label{sec:dlr}

\subsubsection{Railway network in consideration}\label{sec:NetworkDlr}

The second case study uses the Docklands Light Railway (DLR) network shown in Figure~\ref{Figure_loop_docklands}. The DLR is an automated light-metro network in London that operates under the \textsf{SelTrac} CBTC system \citep{selTrac}, and it has served as a benchmark for energy-efficient timetabling. Its topology features loops and turn-around movements, so it tests the model on a structurally different network. The benchmark network contains two directed train-lines, denoted Line~1 and Line~2. Stations are labeled by uppercase abbreviations, and each station has a platform in each direction; for example, BAN has one platform on each line. The labels 1067 and ISP1 denote intermediate stopping points on Line~1, whereas BANH and CROH denote turn-around points on Lines~1 and~2, respectively.

 \begin{figure}[htp]
	\centering
	\psscalebox{0.4 0.4} {
		\begin{pspicture}(0,-8.8386345)(19.471958,8.8386345)
			\psline[linecolor=black, linewidth=0.08, arrowsize=0.05291666666666668cm 2.0,arrowlength=1.4,arrowinset=0.0]{->}(9.479467,-8.288409)(8.219467,-7.5884094)
			\psline[linecolor=black, linewidth=0.08, arrowsize=0.05291666666666668cm 2.0,arrowlength=1.4,arrowinset=0.0]{->}(8.199468,-7.5684094)(8.179467,-6.9684095)
			\psline[linecolor=black, linewidth=0.08, arrowsize=0.05291666666666668cm 2.0,arrowlength=1.4,arrowinset=0.0]{->}(8.159468,-7.0484095)(7.0594673,-6.4084096)
			\psline[linecolor=black, linewidth=0.08, arrowsize=0.05291666666666668cm 2.0,arrowlength=1.4,arrowinset=0.0]{->}(7.0594673,-6.4084096)(7.0594673,-3.5284095)
			\psline[linecolor=black, linewidth=0.08, arrowsize=0.05291666666666668cm 2.0,arrowlength=1.4,arrowinset=0.0]{->}(7.0594673,-3.5284095)(6.699467,-3.0284095)
			\psline[linecolor=black, linewidth=0.08, arrowsize=0.05291666666666668cm 2.0,arrowlength=1.4,arrowinset=0.0]{->}(6.699467,-3.0284095)(0.7194673,-3.0484095)
			\psline[linecolor=black, linewidth=0.08, arrowsize=0.05291666666666668cm 2.0,arrowlength=1.4,arrowinset=0.0]{->}(0.7194673,-3.0284095)(0.079467274,-2.4084096)
			\psline[linecolor=black, linewidth=0.08, arrowsize=0.05291666666666668cm 2.0,arrowlength=1.4,arrowinset=0.0]{->}(0.93946725,7.8515906)(0.93946725,-1.9484096)
			\psline[linecolor=black, linewidth=0.08, arrowsize=0.05291666666666668cm 2.0,arrowlength=1.4,arrowinset=0.0]{->}(0.9194673,-1.8884095)(1.1994673,-2.1884096)
			\psline[linecolor=black, linewidth=0.08, arrowsize=0.05291666666666668cm 2.0,arrowlength=1.4,arrowinset=0.0]{->}(1.1594672,-2.1084096)(8.319468,-2.1484096)
			\psline[linecolor=black, linewidth=0.08, arrowsize=0.05291666666666668cm 2.0,arrowlength=1.4,arrowinset=0.0]{->}(8.299467,-2.1084096)(9.259467,-3.3484097)
			\psline[linecolor=black, linewidth=0.08, arrowsize=0.05291666666666668cm 2.0,arrowlength=1.4,arrowinset=0.0]{->}(9.239468,-3.3084095)(9.259467,-4.6084094)
			\psline[linecolor=black, linewidth=0.08, arrowsize=0.05291666666666668cm 2.0,arrowlength=1.4,arrowinset=0.0]{->}(9.219467,-4.5684094)(10.419467,-5.2484097)
			\psline[linecolor=black, linewidth=0.08, arrowsize=0.05291666666666668cm 2.0,arrowlength=1.4,arrowinset=0.0]{->}(10.379467,-5.22841)(10.399467,-6.3484097)
			\psline[linecolor=black, linewidth=0.08, arrowsize=0.05291666666666668cm 2.0,arrowlength=1.4,arrowinset=0.0]{->}(10.3594675,-6.2884097)(9.219467,-7.0084095)
			\psline[linecolor=black, linewidth=0.08, arrowsize=0.05291666666666668cm 2.0,arrowlength=1.4,arrowinset=0.0]{->}(9.239468,-6.9684095)(9.239468,-7.3084097)
			\psframe[linecolor=black, linewidth=0.08, dimen=outer](9.039467,-3.7284095)(8.759467,-4.2084093)
			\psframe[linecolor=black, linewidth=0.08, dimen=outer](7.5594673,-3.7284095)(7.279467,-4.2084093)
			\psframe[linecolor=black, linewidth=0.08, dimen=outer](0.8394673,0.43159044)(0.55946726,-0.048409548)
			\psframe[linecolor=black, linewidth=0.08, dimen=outer](0.41946727,0.43159044)(0.13946728,-0.048409548)
			\psframe[linecolor=black, linewidth=0.08, dimen=outer](0.85946727,3.5715904)(0.5794673,3.0915904)
			\psframe[linecolor=black, linewidth=0.08, dimen=outer](0.85946727,6.4315906)(0.5794673,5.9515905)
			\psframe[linecolor=black, linewidth=0.08, dimen=outer](0.41946727,6.4315906)(0.13946728,5.9515905)
			\psframe[linecolor=black, linewidth=0.08, dimen=outer](0.3994673,8.131591)(0.11946728,7.6515903)
			\psframe[linecolor=black, linewidth=0.08, dimen=outer](7.5394673,-5.4884095)(7.259467,-5.9684095)
			\psframe[linecolor=black, linewidth=0.08, dimen=outer](10.199468,-5.5084095)(9.919467,-5.9884095)
			\psframe[linecolor=black, linewidth=0.08, dimen=outer](4.339467,-2.6084096)(3.7594674,-2.9084096)
			\psframe[linecolor=black, linewidth=0.08, dimen=outer](4.339467,-2.2084095)(3.7594674,-2.5084095)
			\psframe[linecolor=black, linewidth=0.08, dimen=outer](1.9194673,-2.6484096)(1.3394673,-2.9484096)
			\psframe[linecolor=black, linewidth=0.08, dimen=outer](1.9194673,-2.2084095)(1.3394673,-2.5084095)
			\psframe[linecolor=black, linewidth=0.08, dimen=outer](13.459467,-7.5084095)(12.879467,-7.8084097)
			\psframe[linecolor=black, linewidth=0.08, dimen=outer](10.699468,-7.8884096)(10.119468,-8.18841)
			\psframe[linecolor=black, linewidth=0.08, dimen=outer](16.339468,-7.8684096)(15.759467,-8.168409)
			\psframe[linecolor=black, linewidth=0.08, dimen=outer](16.339468,-7.5084095)(15.759467,-7.8084097)
			\psframe[linecolor=black, linewidth=0.08, dimen=outer](13.459467,-7.8684096)(12.879467,-8.168409)
			\psframe[linecolor=black, linewidth=0.08, dimen=outer](10.699468,-7.5284095)(10.119468,-7.8284097)
			\psframe[linecolor=black, linewidth=0.08, dimen=outer](6.339467,-2.2284095)(5.759467,-2.5284095)
			\psframe[linecolor=black, linewidth=0.08, dimen=outer](18.819468,-7.4684095)(18.239468,-7.7684097)
			\rput[bl](1.0994673,7.8315907){BANH}
			\rput[bl](1.0994673,6.0715904){BAN}
			\rput[bl](1.1594672,3.2115905){1067}
			\rput[bl](1.0594673,0.07159045){SHA}
			\rput[bl](1.3394673,-1.9884095){LIM}
			\rput[bl](3.7194672,-1.9884095){WES}
			\rput[bl](5.699467,-1.9884095){ISP1}
			\rput[bl](9.379467,-4.0084095){WIQ}
			\rput[bl](10.519467,-5.9084096){CAW}
			\rput[bl](10.099467,-8.76841){HEQ}
			\rput[bl](12.839467,-8.708409){SOQ}
			\rput[bl](15.739468,-8.708409){CRO}
			\rput[bl](17.999468,-7.1084094){CROH}
			\rput[bl](10.019467,-7.0284095){\LARGE{Line 1}}
			\rput[bl](5.739467,-8.168409){\LARGE{Line 2}}
			\psline[linecolor=black, linewidth=0.08, arrowsize=0.05291666666666667cm 2.0,arrowlength=1.4,arrowinset=0.0]{->}(0.079467274,-2.3884096)(0.03946728,8.81159)
			\psline[linecolor=black, linewidth=0.08, arrowsize=0.05291666666666667cm 2.0,arrowlength=1.4,arrowinset=0.0]{->}(9.239468,-7.2684097)(10.019467,-7.3484097)
			\psline[linecolor=black, linewidth=0.08, arrowsize=0.05291666666666667cm 2.0,arrowlength=1.4,arrowinset=0.0]{->}(9.999467,-7.3284097)(19.479467,-7.3484097)
			\psline[linecolor=black, linewidth=0.08, arrowsize=0.05291666666666667cm 2.0,arrowlength=1.4,arrowinset=0.0]{->}(19.439466,-7.3284097)(18.779467,-8.248409)
			\psline[linecolor=black, linewidth=0.08, arrowsize=0.05291666666666667cm 2.0,arrowlength=1.4,arrowinset=0.0]{->}(18.759468,-8.248409)(9.459467,-8.30841)
			\psline[linecolor=black, linewidth=0.08, arrowsize=0.05291666666666667cm 2.0,arrowlength=1.4,arrowinset=0.0]{->}(0.03946728,8.81159)(0.93946725,7.8315907)
		\end{pspicture}
	}
	\caption{Docklands Light Railway network used in the numerical study.}
	\label{Figure_loop_docklands}
\end{figure}

\subsubsection{Results of the numerical study}\label{sec:resultsDlr}

We evaluate the model on 17 real-world instances. Each instance corresponds to an operational schedule with a different nominal service headway, ranging from 1~min~22~s to 4~min~17~s; shorter headways allow for more trains, so the instances span 108--336 active trains over a service period of approximately 4.4 hours. The baseline (initial) effective energy consumption of each instance is obtained by evaluating \eqref{eq:predicted_eng} at the existing timetable. As in the Shanghai study, both the baseline and the optimized values are computed with the model, so the comparison uses a single evaluator. The energy coefficients for these instances were fitted network-wide from \textsf{KDC} simulations of representative Docklands Light Railway runs (a single traction fit shared by all trips and one transfer fit per event side), and the robust trip-time windows are 8--44~seconds wide (median 18), with crossing-over windows of 10--11~seconds. Table~\ref{tab:Results-of-the-numerical-study-dlr} summarizes our findings.

\begin{table*}[htbp]
  \centering
  \scriptsize 
  \caption{Summary of numerical results for Docklands Light Railway instances.}
  \label{tab:Results-of-the-numerical-study-dlr}
    \begin{tabular}{C{1.7cm}C{1.7cm}C{2.05cm}C{1.9cm}C{1.9cm}C{1.9cm}C{1.9cm}}
    \toprule
    \textbf{\#Trains} & \textbf{\#Variables} & \textbf{\#Constraints} & \textbf{Solution time (s)} & \textbf{Initial effective energy consumption (kWh)} & \textbf{Predicted final effective energy consumption (kWh)} & \textbf{Reduction predicted by our model (\%)} \\
    \midrule
    108 & 2669 & 10109 & 0.020 & 43317.82 & 34942.62 & 19.33 \\
    121 & 3691 & 12497 & 0.026 & 48029.98 & 37115.47 & 22.72 \\
    144 & 4783 & 15531 & 0.036 & 57145.17 & 43842.99 & 23.28 \\
    156 & 4345 & 15437 & 0.034 & 62554.02 & 49058.68 & 21.57 \\
    180 & 5547 & 18703 & 0.044 & 71421.42 & 55989.84 & 21.61 \\
    192 & 6495 & 20919 & 0.047 & 74739.52 & 57419.50 & 23.17 \\
    204 & 6777 & 22025 & 0.049 & 80164.56 & 61514.73 & 23.26 \\
    216 & 7755 & 24291 & 0.055 & 84901.98 & 63246.11 & 25.51 \\
    228 & 8507 & 26171 & 0.060 & 88791.14 & 67805.80 & 23.63 \\
    240 & 10211 & 29647 & 0.066 & 92279.25 & 69430.24 & 24.76 \\
    252 & 9449 & 29013 & 0.083 & 98841.16 & 73490.23 & 25.65 \\
    263 & 9833 & 30231 & 0.069 & 102097.22 & 77727.41 & 23.87 \\
    275 & 8783 & 29117 & 0.077 & 107981.96 & 83646.10 & 22.54 \\
    288 & 11155 & 33755 & 0.084 & 111463.79 & 82726.07 & 25.78 \\
    300 & 11077 & 34261 & 0.105 & 117099.32 & 89096.60 & 23.91 \\
    311 & 10489 & 33859 & 0.090 & 123214.48 & 95110.66 & 22.81 \\
    336 & 14529 & 41913 & 0.108 & 131531.13 & 93620.88 & 28.82 \\
    \bottomrule 
  \end{tabular}
\end{table*}

\begin{figure}[ht]
    \centering
    \begin{subfigure}[b]{0.4\textwidth}
        \centering
        \includegraphics[width=\textwidth]{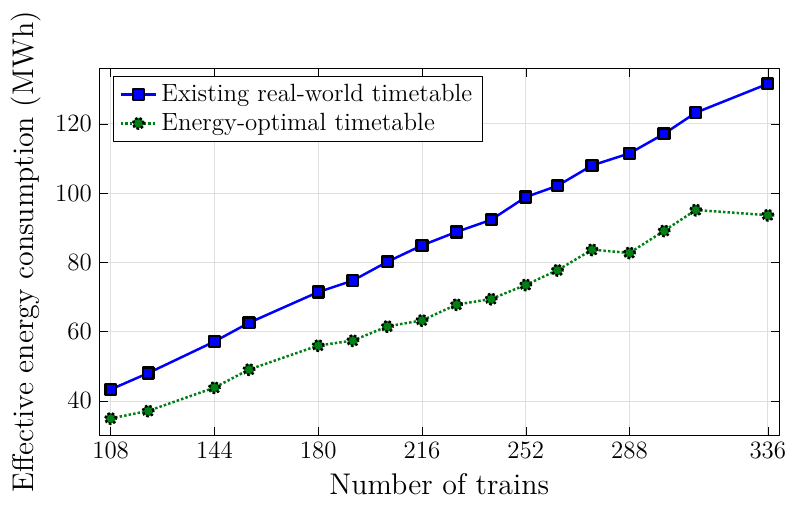}
        \caption{}
        \label{fig:dlr-energy}
    \end{subfigure}
    \hfill
    \begin{subfigure}[b]{0.4\textwidth}
        \centering
        \includegraphics[width=\textwidth]{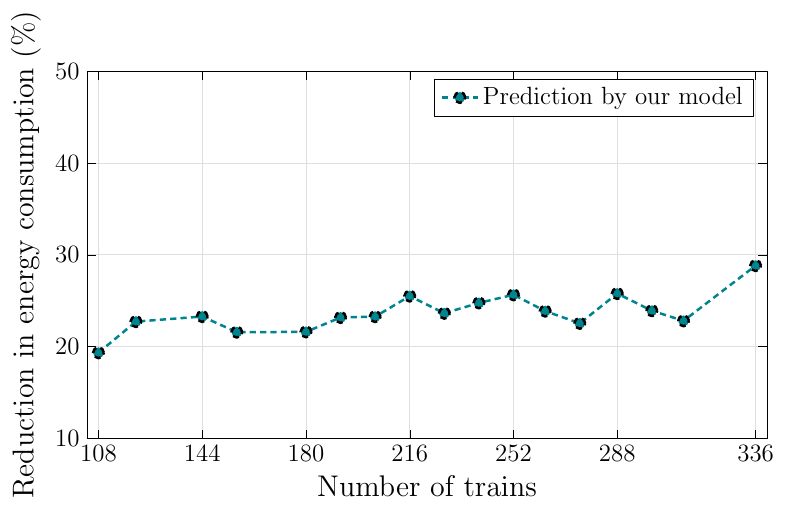}
        \caption{}
        \label{fig:dlr-reduction}
    \end{subfigure}
    \caption{Docklands Light Railway results. (a)~Model-computed effective energy consumption of the existing timetables and of the energy-optimal timetables. (b)~Corresponding percentage reduction predicted by our model.}
    \label{Fig:b_dlr}
\end{figure}
\FloatBarrier

Across the 17 instances, the number of trains ranges from 108 to 336, while the linear programs contain 2,669--14,529 variables and 10,109--41,913 constraints. All minimization solves terminated optimally, and the median solution times range from 0.020 to 0.108 seconds. 

The model-computed energy values in Figure~\ref{fig:dlr-energy} generally increase with the number of trains, but the percentage reductions in Figure~\ref{fig:dlr-reduction} are not monotone because they depend on how energy-efficient the initial timetables are. The predicted savings range from 19.33\% to 28.82\%, indicating that the benefit depends on the synchronization opportunities and timing structure of an instance, not only on its train count. Unlike in the Shanghai study, we did not have sufficient physics-related data required to predict energy savings through \textsf{SPSIM} accurately, so we just present the savings as predicted by our model for this network.

Table~\ref{tab:decomposition-dlr} again decomposes the predicted reductions in the same way as Shanghai study. With sparser traffic the traction component dominates more strongly, supplying 61.7\%--89.8\% of the reduction, yet synchronization still contributes 10.2\%--38.3\%: transferred energy grows by a factor of 2.6--6.8 over the baselines, and in the two sparsest instances (nominal headways 2~min~58~s and 4~min~17~s) the baseline timetables does not allow for any transfer at all, so the entire transferred-energy term is created by our optimization model

\begin{table*}[htbp]
  \centering
  \scriptsize
  \caption{Decomposition of the predicted effective-energy reduction for the Docklands Light Railway instances, in the format of Table~\ref{tab:decomposition-sl8}.}
  \label{tab:decomposition-dlr}
  \begin{tabular}{r rr rr rr}
    \toprule
    & \multicolumn{2}{c}{Baseline} & \multicolumn{2}{c}{Optimized} & \multicolumn{2}{c}{Share of reduction (\%)} \\
    \cmidrule(lr){2-3}\cmidrule(lr){4-5}\cmidrule(lr){6-7}
    \#Trains & $E^{\textrm{con}}$ & $E^{\textrm{reg}}$ & $E^{\textrm{con}}$ & $E^{\textrm{reg}}$ & traction & transfer \\
    \midrule
    108 & 43317.82 & 0.00 & 35797.28 & 854.66 & 89.8 & 10.2 \\
    121 & 48499.51 & 469.54 & 40072.60 & 2957.13 & 77.2 & 22.8 \\
    144 & 57708.61 & 563.44 & 47681.28 & 3838.29 & 75.4 & 24.6 \\
    156 & 62554.02 & 0.00 & 51691.04 & 2632.36 & 80.5 & 19.5 \\
    180 & 72099.41 & 677.99 & 59565.28 & 3575.44 & 81.2 & 18.8 \\
    192 & 76944.82 & 2205.30 & 63575.04 & 6155.55 & 77.2 & 22.8 \\
    204 & 81790.23 & 1625.66 & 67584.97 & 6070.24 & 76.2 & 23.8 \\
    216 & 86635.64 & 1733.66 & 71596.67 & 8350.56 & 69.4 & 30.6 \\
    228 & 91335.62 & 2544.48 & 75459.02 & 7653.22 & 75.7 & 24.3 \\
    240 & 96181.02 & 3901.77 & 79468.79 & 10038.55 & 73.1 & 26.9 \\
    252 & 101026.43 & 2185.27 & 83478.58 & 9988.35 & 69.2 & 30.8 \\
    263 & 105390.13 & 3292.91 & 87077.48 & 9350.07 & 75.1 & 24.9 \\
    275 & 110235.53 & 2253.58 & 91087.55 & 7441.45 & 78.7 & 21.3 \\
    288 & 115417.23 & 3953.44 & 95362.62 & 12636.56 & 69.8 & 30.2 \\
    300 & 120262.64 & 3163.32 & 99380.09 & 10283.49 & 74.6 & 25.4 \\
    311 & 124626.33 & 1411.85 & 102971.55 & 7860.88 & 77.1 & 22.9 \\
    336 & 134653.43 & 3122.30 & 111256.31 & 17635.42 & 61.7 & 38.3 \\
    \bottomrule
  \end{tabular}
\end{table*}

\section{Architectural framework for industrial integration of the proposed model}\label{sec:archframe}

\begin{figure}[htbp]
\begin{centering}
\includegraphics[width=0.94\textwidth]{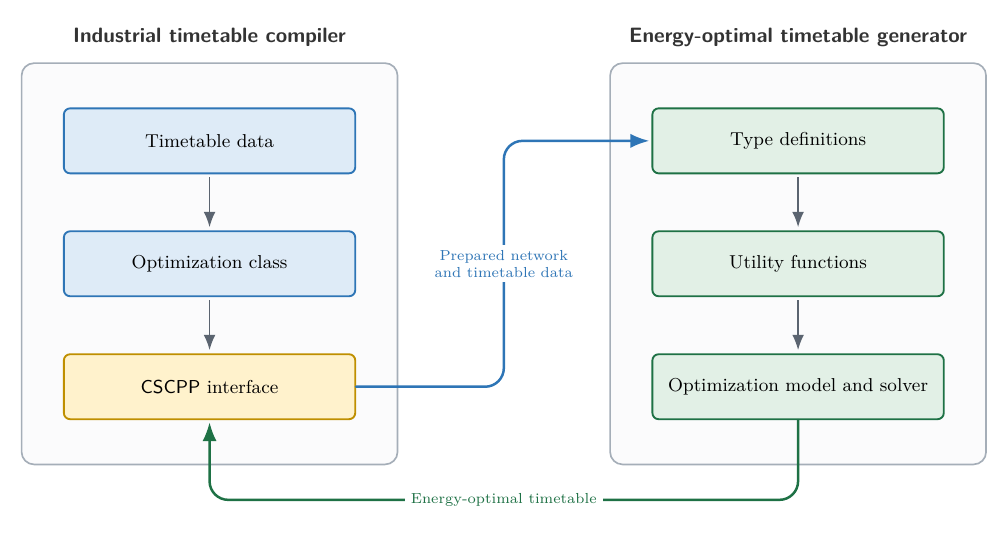}
\par\end{centering}
\caption{Architecture for integrating the energy-optimal timetable generator with an industrial timetable compiler. The \textsf{CSCPP} interface sends prepared network and timetable data to the generator and returns the resulting energy-optimal timetable to the compiler.\label{fig:julia_embedd}}
\end{figure}

This section describes the architectural framework for incorporating the proposed optimization model into Hitachi Rail Canada's industrial timetable compiler. Figure~\ref{fig:julia_embedd} separates the existing timetable-compiler components from the energy-optimal timetable generator and shows the two data flows between them.

\FloatBarrier

For the purpose of this integration framework, the existing timetable compiler is partitioned into three discrete modules:

\begin{itemize}
\item \emph{Timetable data:} This module houses the necessary data to solve the optimization problem. Additionally, it includes other information specific to the railway network under consideration, such as gradient attributes and speed limitations in various track segments, all presented in raw form.

\item \emph{Optimization class:} This module transforms the raw data into a structure that is compatible with the optimization algorithm we employ.

\item \emph{\textsf{CSCPP} interface:} This interface is the compiler's data-exchange layer; it provides bidirectional communication between the industrial timetable compiler and the energy-optimal timetable generator. It sends the prepared network and timetable data to the generator, retrieves the energy-optimal timetable, and converts the result into a format suitable for deployment in a CBTC system.
\end{itemize}

The energy-optimal timetable generator is similarly sub-divided into three modules:

\begin{itemize}
\item \emph{Type definitions:} This module contains the data structures essential for describing a railway network, as described in Sections \ref{sec:Modelling-the-constraints} and \ref{sec:Final-optimization-model}.
 
\item \emph{Utility functions:} This module constructs the optimization model \eqref{eq:final_lp_model}, predicts effective energy consumption via \eqref{eq:predicted_eng}, and transforms the output timetable into a \textsf{CSCPP}-compatible representation.

\item \emph{Optimization model and solver:} This module contains the linear program to solve \eqref{eq:final_lp_model}. The resulting timetable is returned to the \textsf{CSCPP} interface through the utility functions.
\end{itemize}

\section{Conclusion\label{Conclusion} }

We presented a single-stage, data-driven linear programming model that designs energy-optimal timetables for CBTC-operated metro railway networks by minimizing effective energy consumption: the total traction energy minus the regenerative braking energy transferred from braking trains to accelerating trains. The model embeds the transfer of regenerative energy in linear hypograph constraints and predicts the effective energy consumption of a timetable without physics-based simulation. On 11 full-day operational instances of Shanghai Metro Line 8 with up to 1,332 active trains, it computed energy-optimal timetables in under 0.7 seconds per instance and predicted effective-energy reductions of 19.52\%--28.17\% relative to the existing operational timetables, in close agreement with the physics-based simulator \textsf{SPSIM}. On 17 Docklands Light Railway real-world instances with up to 336 active trains, it solved in under 0.11 seconds per instance and predicted reductions of 19.33\%--28.82\%.

Our model has some limitations. The affine models of traction and regenerative energy are justified only over the narrow trip-time windows of Assumption~\ref{assum:triptime} that is typical of CBTC systems; however substantially wider windows that is common in inter-city travels may require nonlinear fits. The FWHM rectangles approximate the power-versus-time profiles conservatively, so the model tends to underestimate transferred regenerative energy relative to \mbox{\textsf{SPSIM}}'s exact overlap computation. Also, the Docklands Light Railway predictions could not be validated against a physics-based simulator such as \textsf{SPSIM} due to lack of appropriate data required for the simulator for Docklands Light Railway. 

The model is being integrated into Hitachi Rail Canada's industrial timetable compiler through the architecture described in Section~\ref{sec:archframe}. Natural directions for future work include extending the model to handle stochastic disturbances and delay recovery, and incorporating passenger-demand information into the timetable design.

\bibliographystyle{plainnat}
\bibliography{train_opt} \end{document}